\documentstyle{amlts}
\begin{document}
\annalsline{156}{2002}
\received{January 19, 2001}
\startingpage{295}
\def\bye{\end{document}}
 \font\tenrm=cmr10
\input amssym.def
\input amssym.tex
\def\ritem#1{\item[{\rm #1}]}
\input boxedeps.tex 
\SetepsfEPSFSpecial 
\HideDisplacementBoxes
\def\figin#1#2{
$$
 {\BoxedEPSF{#1  scaled
#2}%
}%
$$
\noindent}
\def\rmg{\hbox{$\frak g\hskip-6pt {\frak g}$}\hskip1pt} 
\def\rmsl{\hbox{$\frak s\hskip-5.1pt \frak s \hskip-5.15pt \frak s \hskip.5pt
\frak l\hskip-3.25pt\frak l\hskip-2.75pt \frak l\hskip-3pt \lower.25pt\hbox{$\frak l$}$}\hskip1pt}
\def\rmu{\hbox{$\frak u\hskip-6pt\frak u$}\hskip1pt}
\def\rmm{\hbox{$\frak m\hskip-8pt\frak m$}\hskip1pt}
\catcode`\@=11
\font\twelvemsb=msbm10 scaled 1100
\font\tenmsb=msbm10
\font\ninemsb=msbm10 scaled 800
\newfam\msbfam
\textfont\msbfam=\twelvemsb  \scriptfont\msbfam=\ninemsb
  \scriptscriptfont\msbfam=\ninemsb
\def\msb@{\hexnumber@\msbfam}
\def\Bbb{\relax\ifmmode\let\next\Bbb@\else
 \def\next{\errmessage{Use \string\Bbb\space only in math
mode}}\fi\next}
\def\Bbb@#1{{\Bbb@@{#1}}}
\def\Bbb@@#1{\fam\msbfam#1}
\catcode`\@=12

 \catcode`\@=11
\font\twelveeuf=eufm10 scaled 1100
\font\teneuf=eufm10
\font\nineeuf=eufm7 scaled 1100
\newfam\euffam
\textfont\euffam=\twelveeuf  \scriptfont\euffam=\teneuf
  \scriptscriptfont\euffam=\nineeuf
\def\euf@{\hexnumber@\euffam}
\def\frak{\relax\ifmmode\let\next\frak@\else
 \def\next{\errmessage{Use \string\frak\space only in math
mode}}\fi\next}
\def\frak@#1{{\frak@@{#1}}}
\def\frak@@#1{\fam\euffam#1}
\catcode`\@=12


 
\title{Parametrizing nilpotent
orbits\\ via Bruhat-Tits theory} 
\shorttitle{Parametrizing nilpotent orbits} 

 \acknowledgements{Supported by National Science
Foundation Postdoctoral Fellowship 98-04375.\hfill\break
\phantom{hihow}1991 {\it Mathematics Subject Classification}. Primary 20G25; Secondary 17B45, 20G15.\hfill\break
\phantom{hihow}{\it Key words and phrases}.  Bruhat-Tits building, nilpotent orbit, reductive group.}
\author{Stephen DeBacker}
\institutions{The University of Chicago, Chicago, IL\\
{\eightpoint {\it Current address\/}:} Harvard University, Cambridge, MA\\
{\eightpoint {\it E-mail address\/}: debacker@math.harvard.edu
}}

\newcommand{\inv}{^{-1}}
\let\without\smallsetminus		
\newcommand{\X}{{\protect{\bf X}}}

\newcommand{\Z}{{{\Bbb Z}}}
\newcommand{\R}{{{\Bbb R}}}
\newcommand{\N}{{{\Bbb N}}}
\newcommand{\Q}{{{\Bbb Q}}}
\newcommand{\F}{{{\Bbb F}}}

\newcommand{\supp}{\mathop{\rm  supp}\nolimits}
\newcommand{\stab}{\mathop{\rm stab}\nolimits}
\newcommand{\End}{\mathop{\rm End}\nolimits}
\newcommand{\SL}{\mathop{\rm SL}\nolimits}
\newcommand{\Gal}{\mathop{\rm Gal}\nolimits}
\newcommand{\Lie}{\mathop{\rm Lie}\nolimits}

\newcommand{\glo}{{{\rm GL}_1}}
\newcommand{\sllie}{{\mathbold{{\frak   sl}}}}

\newcommand{\Ind}[3]{\operatorname{Ind}_{#1}^{#2}{#3}}

\newcommand{\pmat}[4]{%
	\left(\begin{matrix}
	\vphantom{P^2}#1 & \vphantom{P^2}#2 \\
	\vphantom{P^2}#3 & \vphantom{P^2}#4
	\end{matrix}\right)}

\let\oldwp\wp
\def\wp{{\raise.4ex\hbox{\hbox{$\oldwp$}}}}

\let\Ind\relax
\newcommand{\Ind}{\mathop{\rm Ind}\nolimits}
\newcommand{\dAd}{\mathop{\rm Ad}\nolimits}
\newcommand{\dad}{\mathop{\rm ad}\nolimits}
\newcommand{\ad}{\mathop{\rm ad}\nolimits}
\newcommand{\aut}{\mathop{\rm Aut}\nolimits}
\newcommand{\dist}{\mathop{\rm dist}\nolimits}
\newcommand{\dorbit}{{{\cal  O}}}
\newcommand{\dborbit}{{{\cal  O}}}
\newcommand{\zz}{{{\frak   z}}}
\newcommand{\nn}{{{\frak   n}}}
\newcommand{\pp}{{{\frak   p}}}
\newcommand{\uu}{{{\frak   u}}}
\renewcommand{\gg}{{\frak   g}} 
\newcommand{\bgg}{{\mathbold{\gg}}}
\newcommand{\buu}{{\mathbold{\uu}}}
\newcommand{\mm}{{{\frak   m}}}
\renewcommand{\AA}{{{\cal  A}}}
\newcommand{\BB}{{{\cal  B}}}
\newcommand{\NN}{{{\cal  N}}}
\newcommand{\UU}{{{\cal  U}}}

\newcommand{\bG}{{\bf G}}
\newcommand{\bGn}{{\bf G^\circ}}
\newcommand{\Gn}{{G^\circ}}
\newcommand{\bS}{{\bf S}}
\newcommand{\bT}{{\bf T}}
\newcommand{\bZ}{{\bf Z}}
\newcommand{\bM}{{\bf M}}
\newcommand{\bmm}{\mathbold{{\frak   m}}}

\newcommand{\ff}{{\frak   f}}
\newcommand{\ffc}{{\frak   F}}

\newcommand{\bfZ}{{\bf \mathsf{Z}}}
\newcommand{\bfG}{{\bf \mathsf{G}}}
\newcommand{\bfT}{{\bf \mathsf{T}}}
\newcommand{\bfB}{{\bf \mathsf{B}}}
\newcommand{\bfP}{{\bf \mathsf{P}}}
\newcommand{\bfS}{{\bf \mathsf{S}}}

\newcommand{\lsup}[1]{{}^{#1}}
\newcommand{\dabs}[1]{\left|#1\right|}

\newcommand{\ident}{\mathop=\limits^{i}}

 \centerline{\bf Abstract}
\vglue9pt

Let $k$ denote a field with nontrivial discrete valuation.   We assume that $k$
is complete with perfect residue field.   Let $G$ be the
group of $k$-rational points of a reductive, linear algebraic group
defined over $k$. Let $\gg$ denote 
the Lie algebra of $G$. Fix $r \in \R$.   Subject to some
restrictions, we show that the set of distinguished degenerate 
Moy-Prasad cosets of depth $r$ (up to an equivalence relation)
parametrizes the 
nilpotent orbits in $\gg$. 
 \vglue-12pt

\section{Introduction}

In this paper we give a uniform
parametrization of   the nilpotent orbits in the Lie algebra of a
$p$-adic reductive group.
This classification, which was motivated by harmonic analysis
considerations,  matches nilpotent orbits with certain
equivalence classes that arise naturally from Bruhat-Tits
theory. 

\demo{{\rm 1.1.} Motivation}
In the early 1970s Harish-Chandra and Roger Howe studied the local
behavior of the character of an irreducible smooth representation of
a reductive $p$-adic group~\cite{hc:queens}, \cite{howe:fourier}.  For
example, they established what is now called the Harish-Chandra-Howe
local character expansion -- in some unspecified neighborhood of the
identity the character can be expressed as a linear combination of the
Fourier transforms of nilpotent orbital integrals.  At the heart of
their proofs was a remarkable finiteness statement, referred to as
``Howe's conjecture''~\cite{howe:twoconj}, about invariant
distributions on the Lie algebra.
In some stunning work of the 1990s, J.-L.~Waldspurger proved a very
precise version of Howe's conjecture for ``unramified 
classical groups''~\cite{waldspurger:finitude}.  This sharpened
finiteness statement allowed him to relate the range of validity for
the Harish-Chandra-Howe local character expansion to the first
occurrence of fixed-vectors with respect to congruence filtration
subgroups~\cite{waldspurger:homogeneity}. 

The fundamental work of Allen~Moy and
Gopal~Prasad~\cite{moy-prasad:jacquet}, \cite{moy-prasad:K-types} introduced
new ways to use
the structure theory of F.~Bruhat and
J.~Tits~\cite{bruhat-tits:one}, \cite{bruhat-tits:two} to study questions in
representation  theory.  One consequence of their work is that to each
representation we can attach a number, called the  depth of the
representation. Roughly speaking, this number measures the first occurrence
of fixed-vectors with respect to all the natural subgroup filtrations
arising from Bruhat-Tits theory.   The conjecture of Thomas Hales,
Allen Moy, and Gopal Prasad~\cite[\S1]{moy-prasad:K-types} seeks to strengthen
the results of J.-L.~Waldspurger by asking
if the range of validity for the Harish-Chandra-Howe local character
expansion is controlled by the depth of the  
representation; such a result
would greatly enhance our understanding of characters. The
parametrization of nilpotent orbits presented in this article is
the cornerstone  of my proof of their conjecture. The remainder of the
proof appears
in~\cite{adler-debacker:kirilovn}, \cite{debacker:aproofn}, \cite{debacker:Gp}. 
\enddemo

 1.2. {\it The parametrization}.
In a special situation ($r = 0$), the
main result of this 
paper may be viewed as an affine
analogue of Bala-Carter
theory [4], [5].  Namely,
it provides a classification of 
the nilpotent orbits  
in terms of equivalence classes of pairs $(G_F/G_F^+,X)$.  Here $F$ is
a facet in the Bruhat-Tits building of our group, $G_F$ is
the associated parahoric subgroup with pro-unipotent radical
$G_F^+$, and $X$ 
is a distinguished element of the Lie algebra of $G_F/G_F^+$.  (Recall
that $X$ is called distinguished provided that it is 
nilpotent and does not
lie in a proper Levi subalgebra.)

In this article we prove this special case ($r=0$) and take it one
step further -- we
classify the nilpotent orbits in terms of Moy-Prasad cosets of 
an arbitrary fixed depth $r$ (see below). 
We now  discuss the parametrization scheme in detail.

Let $k$ denote a field with nontrivial discrete valuation.   We assume that $k$
is complete with perfect residue field $\ff$.   
Let $G$ denote the group of $k$-rational
points of a reductive, linear algebraic group $\bG$ defined over
$k$ and let $\gg$ 
denote its Lie algebra.  We let $\Gn$ denote the group of $k$-rational
points of the identity component
$\bGn$ of $\bG$.  Let $\BB(G)$ denote the Bruhat-Tits building
of $\Gn$.  For each pair $(x,r) \in \BB(G) \times \R$, Allen~Moy and
Gopal~Prasad~\cite{moy-prasad:jacquet}, \cite{moy-prasad:K-types} have defined
the (Moy-Prasad) lattices $\gg_{x,r^+} 
\subset \gg_{x,r}$ of $\gg$.  For $x \in \BB(G)$, an element of
$\gg_{x,r}/\gg_{x,r^+}$ is called a  Moy-Prasad coset of depth $r$.

Suppose $r \in \R$.  We partition $\BB(G)$ into generalized
$r$-facets --  two points $x$ and $y$ in $\BB(G)$ belong to the same
generalized 
$r$-facet provided that $\gg_{x,r} = \gg_{y,r}$ and $\gg_{x,r^+} =
\gg_{y,r^+}$.  If $F^*$ is a generalized $r$-facet and $x \in F^*$,
then we define 
the  $\ff$-vector space $V_{F^*} = \gg_{x,r}/\gg_{x,r^+}$.  For
example, if $r = 
0$, then generalized $0$-facets are facets in the usual sense, and if
$F$ is a facet of $\BB(G)$, then $V_{F}$ is $\Lie(G_F/G_F^+)$.

Let $I_r$ denote the set of pairs $(F^*,v)$ where $F^*$ is a
generalized $r$-facet and 
$v$ is an element of $V_{F^*}$.  The set $I_r$  parametrizes the
set of Moy-Prasad cosets of depth $r$.
In~Section~3.6 we define on $I_r$ an equivalence relation,
denoted $\sim$,  which is a  
natural extension of the concept of
associate~\cite{moy-prasad:jacquet}, \cite{moy-prasad:K-types}.

 A pair $(F^*,v) \in I_r$ is degenerate if the
coset it parametrizes contains a nilpotent element.  Let $I_r^n$ denote the subset of $I_r$ consisting of degenerate pairs. 
With some restrictions on $k$ and $\bG$ (see~Section~4.2), we
generalize a result of 
Dan~Barbasch and Allen~Moy~\cite[\protect{\S3}]{barbasch-moy:local}. We
show that 
to each element $(F^*,e)$ of $I_r^n$ we can associate a unique
nilpotent orbit $\dorbit(F^*,e)$.
This orbit is characterized by the fact that it is the
nilpotent orbit of minimal dimension having nontrivial intersection
with the coset
corresponding to $(F^*,e)$.

The set $I_r^n$ is too large for our purposes.  We therefore restrict
our attention to the subset $I_r^d$ of distinguished elements of
$I_r^n$ (see~Section~5.5).  For example, if $r = 0$, then
$(F,e) \in I^n_0$ is 
distinguished if $e$ is a distinguished element of $V_F= \Lie(G_F/G_F^+)$
in the sense discussed above.

We now state Theorem~\ref{thm:parametrize}, the main result of
this paper.  Let $\dorbit (0)$ 
denote the set of nilpotent orbits in $\gg$.
\enddemo

{\elevensc Theorem}.
{\it Assume that all of the hypotheses of Section~{\rm 4.2} hold{\rm .}
There is a bijective correspondence between $I_r^d / \!\sim$ and
$\dorbit(0)$ given by the map which sends $(F^*,e)$ to} $\dorbit(F^*,e)${\rm .}
\vglue4pt

We remark that this result is false without some restrictions on $k$
and $\bG$.  For example, if $k$ is the field of Laurent series over
the field with two elements, then for the group
${\bf SL}_2(k)$ the set $I_0^d / \!\sim$ has cardinality three,
but $\dorbit(0)$ has infinitely many elements.  On the other hand, if
we are not interested in a proof which works in a general setting, then we can
get by with less severe restrictions. For example, we expect that the
theorem is 
true for 
${\bf GL}_n(k)$ with no restrictions on $k$; if $r = 0$, then this
is easy to verify.  If we assume that the residual
characteristic of $k$ is not two, then we expect that the result is
valid for 
split classical 
groups.

In the special case when $r=0$, the parametrization scheme discussed
in this article is inherent (though neither stated nor proved) in
a  paper of\break Dan~Barbasch and Allen~Moy~\cite{barbasch-moy:local}.
Magdy~Assem pointed this out to Robert Kottwitz who, in turn, pointed
it out to me.  Also in the case when $r=0$,\break
J.-L.~Waldspurger~\cite{waldspurger:quelque}  develops a
conjectural parametrization scheme similar to that given here but for
unipotent orbits.  He verifies his conjecture in a
number of cases.  Finally, if $r=0$, $k$ is the field of Laurent series over the complex numbers, and ${\bf G}$ is a connected,
simple, adjoint, and $k$-split group, then the main result of Eric Sommers' paper \cite{sommers:thesis} is equivalent to the
main result of this paper; the proofs, however, are very different. 

I thank both Robert Kottwitz and Gopal Prasad for their many
corrections and improvements to earlier versions of this paper. 
I thank Eugene Kushnirsky and Gopal~Prasad for allowing me to use
their proofs (Lemma~\ref{lem:kushnirsky} and Lemma~\ref{lem:mpalmost},
respectively).
This paper has benefitted from discussions with Jeff~Adler, Robert~Kottwitz,
Allen~Moy, Fiona~Murnaghan, Amritanshu~Prasad, Gopal~Prasad,
Paul~J.~Sally,~Jr., and
Jiu-Kang~Yu.  It is a true pleasure to thank all of these people.  

\vglue-18pt

\section{Notation}
 
2.1. {\it Basic notation}.  
Let $k$ denote a field with nontrivial discrete valuation~$\nu$.    We also denote by $\nu$ the unique
extension of $\nu$ to any algebraic extension of $k$.
We assume
that $k$ is complete and the residue field $\ff$ is 
perfect. Denote the ring of integers of $k$ by $R$ and fix a
uniformizer $\varpi$.

Let $K$ be a fixed maximal unramified extension of $k$.  Let $R_K$
denote the ring of integers of $K$
and let $\ffc$ denote the residue
field of $K$.  Note that $\ffc$ is an algebraic closure of $\ff$.

If $\ff$ has 
positive characteristic, then we let $p$ denote the characteristic of~$\ff$.  If $\ff$ has characteristic zero, then we let $p =
\infty$. Suppose $n \in \Z$.  If  $p <
\infty$, then $(n,p) = \gcd(n,p)$.  If $p = \infty$, then $(n,p) = 1$. 

Let $\bG$ be a reductive,
linear algebraic group defined over $k$. Let $\bGn$ denote the
identity component of $\bG$.  Note the $\bGn$ is a connected,
reductive, linear algebraic group which is defined over $k$.  We let $G =
\bG(k)$ and $\Gn = \bGn(k)$.  We denote by ${\rmg}$ 
  the {L}ie
algebra of $\bG$.   We let $\gg = {\rmg }(k)$, the vector space of
$k$-rational points of ${\rmg }$.  Let $(\!X,Y\!)\!  \mapsto \! [X,Y]$ denote the
Lie algebra product for~$\gg$.

We adopt the following conventions.  We call a subgroup of $\bG$ a
{\it parabolic subgroup of $\bG$} provided that it is a parabolic
subgroup of $\bGn$.  Similar notation applies to tori and Levi subgroups.

Let $L$ be the minimal Galois extension of $K$ such that $\bGn$ is $L$-split.
As in~\cite{moy-prasad:jacquet}, we define $\ell = [L:K]$, and we
normalize $\nu$ by requiring $\nu (L^\times) = \Z$.

If $g \in G$ and $X \in \gg$, then $\lsup{g}X = \dAd(g)X$.  If $X \in
\gg$, then $\lsup{G}X$ denotes the $G$-orbit of $X$ in $\gg$.  We let
$\X_*^k(\bG)$ denote the set of one-parameter $k$-subgroups of $\bG$.

An element $X \in \gg$ is ${\it nilpotent}$ if and only if there
exists $\lambda \in \X_*^k(\bG)$ such that $\lim_{t \rightarrow 0}
\lsup{\lambda(t)}X = 0$.  Let $\NN$ denote the set of nilpotent
elements in $\gg$ and let $\dborbit(0)$ denote the set of
nilpotent $G$-orbits in $\gg$.  
It is more usual to say that an element is nilpotent
if the Zariski closure of its $G$-orbit
contains zero.  Let $\NN''$ denote the set of elements in $\gg$
that are nilpotent in this sense.  We will let $\NN'$ denote the set
of 
elements in $\gg$ which contain zero in the $p$-adic closure of their
$G$-orbit.
It follows that $\NN\subseteq \NN'\subseteq \NN''$.
From~\cite{kempf:instability} we have $\NN = \NN''$ if $k$ is perfect.
From~\cite{adler-debacker:kirilovn} we have that if $k$ is
perfect or 
$\ff$ is finite, then
$\NN$ = $\NN'$.

Similarly, we say that $h \in G$ is
${\it unipotent}$ provided that there 
exists $\lambda \in \X_*^k(\bG)$ such that $\lim_{t \rightarrow 0}
\lambda(t) h (\lambda  (t))\inv = 1$.  

As in~\cite[\S2.2.1]{corvallis} a subset $H$ of $G$ is {\it bounded}
provided that for every $k$-regular function $f$ on $G$, the set
$\nu\bigl(f(H)\bigr)$ is bounded from below.

\demo{{\rm 2.2.} Apartments{\rm ,} buildings{\rm ,} and associated notation} 
Let $\BB(G)$ denote the (enlarged) Bruhat-Tits building of $\Gn$;
 {i.e.}, $\BB(G)$ takes into account the center of $\Gn$.
We identify $\BB(G)$ with the $\Gal(K/k)$-fixed points of
$\BB(\bG,K)$, the Bruhat-Tits building of $\bGn(K)$.

For $\Omega \subset \BB(G)$, we let $\stab_G(\Omega)$ denote the
stabilizer of $\Omega$ in $G$.

We let $\dist \colon \BB(G) \times \BB(G) \rightarrow \R_+$ denote
a (nontrivial) $G$-invariant distance function as discussed
in~\cite[\S2.3]{corvallis}.
For $x, y \in \BB(G)$, let $[x,y]$ denote the geodesic in $\BB(G)$
from $x$ to $y$ and let $(x,y]$ denote $[x,y] \smallsetminus \{x\}$.

For a $k$-Levi subgroup $\bM$ of $\bG$, we identify $\BB(\bM,k)$ in
$\BB(\bG,k)$.  There is not a canonical way to do this, but every
natural embedding of $\BB(\bM,k)$ in $\BB(\bG,k)$ has the same image.

Given a maximal 
$k$-split torus $\bS$  of $\bG$ we have the torus $S = \bS(k)$
in $G$ and the
corresponding apartment $\AA(S)
= \AA(\bS, 
k)$ in $\BB(G)$.  For $\Omega \subset \AA(S)$, we let $A\bigl(\Omega, \AA(S)\bigr)$ denote the
smallest 
affine subspace of $\AA(S)$ containing $\Omega$. 

We let $\Phi (S) = \Phi(\AA) = \Phi (\bS, k)$ denote the set of roots
of $\bG$ with respect to $k$ and $\bS$; we denote by $\Psi (S) = \Psi
(\AA) = \Psi(\bS, k, \nu)$ the set of affine roots of $\bG$ with
respect to $k$, $\bS$,  and $\nu$.  If $\psi \in \Psi(\AA)$,  
then $\dot{\psi} \in \Phi(\AA)$ denotes the gradient of $\psi$.

For $\psi \in \Psi(\AA)$, let $U_\psi$ and $U_{\psi}^+ := U_{\psi^+}$
denote the 
corresponding subgroups of the root group $U_{\dot{\psi}}$
(see~\cite[\S2.4 and \S3.1]{moy-prasad:K-types}).

For $x \in \BB(G)$, we will denote the parahoric subgroup of $\Gn$ attached to
$x$ by 
$G_x$, and we denote its 
pro-unipotent radical by $G_x^+$.  Note that both $G_x$ and
$G_x^+$ depend only on the facet of $\BB(G)$ to which $x$ belongs.
If $F$ is a facet in $\BB(G)$ and $x \in F$, then 
we define $G_F = G_x$ and $G_F^+ = G_x^+$.

Suppose $x \in \BB(G)$.
The quotient $G_x / G_x^+$ is the group of $\ff$-rational points of a
connected reductive group $\bfG_x$ defined over $\ff$.  
We let $\bfZ_x$ denote the $\ff$-split torus in the center of $\bfG_x$
corresponding to the maximal $k$-split torus in the center of
$\bG$.

  We denote the parahoric subgroup of $\bGn(K)$ corresponding to $x \in
   \BB(\bG,K)$ by $\bG(K)_x$.  We denote the pro-unipotent radical of
   $\bG(K)_x$ by $\bG(K)_x^+$.  The  subgroups $\bG(K)_x$ and
   $\bG(K)_x^+$ depend
   only on the facet of $\BB(\bG,K)$ to which $x$ belongs.  If $F$
   is a facet in $\BB(\bG,K)$ and $x \in F$, then we define $\bG(K)_F
   = \bG(K)_x$ and $\bG(K)_F^+ = \bG(K)_x^+$.
For a facet $F$ in $\BB(\bG,K)$, the quotient $\bG(K)_F/\bG(K)_F^+$ is
   the group of $\ffc$-rational points of a connected, reductive $\ffc$-group
   $\bfG_F$.  
\enddemo

2.3. {\it The Moy-Prasad filtrations of $\gg$}.  When $\ff$ is finite,
in~\cite{moy-prasad:jacquet}, \cite{moy-prasad:K-types} Allen Moy and Gopal Prasad associate to a pair $(x,r) \in
\BB(G)
\times \R$ a lattice $\gg_{x,r}$ in $\gg$.  There is no difficulty in
extending their definition to our setting
(see~\cite{adler-debacker:kirilovn}), and we will not repeat the
definition here. However, we
will need to know that $\gg_{x,r}$ has a 
nice decomposition (with respect to the field $k$).

Suppose that $\bS$ is a maximal $k$-split torus of $\bG$.  Let $\bT$ be a
maximal $K$-split $k$-torus containing $\bS$.  
We identify $\AA(\bS,k)$ with $\AA(\bT,K)^{\Gal(K/k)}$.  For $\phi \in
\Psi\bigl(\AA(\bT,K)\bigr)$, we
define as in~\cite[\S3.2]{moy-prasad:K-types} the lattice ${\rmu}_\phi$  
in the root space ${\rmu}_{\dot{\phi}}$ of ${\rmg }(K)$.  For $\psi \in
\Psi\bigl(\AA(\bS,k)\bigr)$, define the lattice $\gg_\psi$ in the root space
$\gg_{\dot{\psi}}$ of $\gg$ to be the $\Gal(K/k)$-fixed points of 
$$\bigoplus_{\phi \in \Psi(\AA(\bT,K)); \, \, \, \phi |_{\AA(\bS,k)} =
\psi} {\rmu}_\phi.$$ 
One can check that for $\psi, \psi' \in \Psi\bigl(\AA(\bS,k)\bigr)$ we have $\gg_\psi
= \gg_{\psi'}$ if and only if $\psi = \psi'$.  We also define the
lattice $\gg_\psi^+$ in the root space
$\gg_{\dot{\psi}}$ by
$$\gg_\psi^+ = \bigcup \gg_{\psi'},$$
where the union is over those affine
roots $\psi' \in \Psi \bigl( \AA(\bS,k) \bigr)$ such that $\dot{\psi'}
= \dot{\psi}$ and $\psi'(x) > \psi(x)$ for some (hence any) $x \in
\AA(\bS,k)$.
\pagegoal=48pc

Let ${\rmm}$ denote the {L}ie algebra of the $k$-Levi subgroup
$C_{\bGn}(\bS)$.  Let $\mm = 
{\rmm}(k)$.
For $x \in \AA(\bS,k)$, let $\mm_r = \mm \cap
\gg_{x,r}$.  The lattice $\mm_r \subset \mm$ is independent of the choice of $x \in
\AA(\bS,k)$.  If $x \in \AA(\bS,k)$, then
$$\gg_{x,r} = \mm_r \oplus \sum_{\psi \in \Psi(\AA(\bS,k)); \, \, \, \psi(x)
\geq r} \gg_\psi.$$
We define $\gg_{x,r^+} := \cup_{s > r} \gg_{x,s}$.

For $x \in \BB(\bG,K)$ and $s \in \R$, we denote by ${\rmg }(K)_{x,s}$
the Moy-Prasad filtration lattice of ${\rmg }(K)$ associated to $x$ and
$s$.  If $x$ is $\Gal(K/k)$-invariant, then $\gg_{x,s} =
(\hskip1pt{\rmg }(K)_{x,s})^{\Gal(K/k)}$. 

For $(x,r) \in \BB(G) \times \R_{\geq0}$, Moy and Prasad also define
subgroups $G_{x,r} \subset G_x$ (see
also~\cite{adler-debacker:kirilovn}).

\section{Generalized $r$-facets and associated objects} \label{sec:generalized}
 
Fix $r \in \R$.  None of the statements in this section depend on the
structure of $\gg$ as a Lie algebra.  Consequently, all statements
remain true when the roles of $\gg$ and $\gg^*$ are interchanged.

\demo{{\rm 3.1.} $r$-facets}
  Fix a maximal $k$-split torus $\bS$ of $\bG$.  Let $\AA =
\AA(\bS,k)$ be the corresponding apartment in $\BB(G)$.  For
each $\psi \in \Psi(\AA)$, let  
$$H_{\psi -r} := \{x \in \AA \, | \, \psi(x) = r \}.$$ 
This defines a facet structure on $\AA$;  
a nonempty subset $F_\AA \subset \AA$ is called an {\it $r$-facet of $\AA$}
provided that there exists a finite subset $S \subset \Psi(\AA)$ such
that 
$$F_\AA \subset H_S := \bigcap_{\psi \in S} H_{\psi -r},$$
and $F_\AA$ is a connected component (in $H_S$) of
$$H_S \smallsetminus \bigcup_{\psi \in \Psi(\AA) \smallsetminus S} (H_S
\cap H_{\psi-r}).$$
If $F_\AA$ is an $r$-facet of $\AA$, then we
define the dimension of $F_\AA$ by  
$$\dim F_\AA := \dim A(F_\AA, \AA).$$ 
If $F_\AA$ is an $r$-facet of $\AA$ of maximal dimension, then $F_\AA$
is called an {\it $r$-alcove} of $\AA$.   
 
\demo{Example {\rm 3.1.1}} 
$$\{\hbox{$r$-alcoves of $\AA$} \} = \{ \hbox{connected components of
$\AA \smallsetminus \bigcup_{\psi \in \Psi(\AA)} H_{\psi -r}$}\}.$$ 
\enddemo

{\it Remark} 3.1.2.  
$F_\AA$ is an $r$-facet of $\AA$ if and only if $F_\AA$ is a
$(-r)$-facet of $\AA$. 
\vglue9pt

If $F_\AA$ is an $r$-facet of $\AA$ and $x,y \in F_\AA$, then
$\gg_{x,r} = \gg_{y,r}$ and $\gg_{x,r^+} = \gg_{y,r^+}$.
Therefore,
the following definitions make sense. 
 
\specialnumber{3.1.3}\numbereddemo{Definition}Let $F_\AA$ be an $r$-facet of $\AA$.  Fix $x \in F_{\AA}$. 
$$\gg_{F_\AA} := \gg_{x,r}$$ 
and 
$$\gg_{F_\AA}^+ := \gg_{x,r^+}.$$ 
\enddemo 
 
Sometimes, in order to avoid confusion, we denote $\gg_{F_\AA}$ by
$\gg_{F_\AA,r}$ and 
$\gg_{F_\AA}^+$ by $\gg_{F_\AA, r^+}$.

\specialnumber{3.1.4}\proclaim{Lemma} \label{lem:rfacec} 
Let $F_\AA$ be an $r$\/{\rm -}\/facet of $\AA${\rm .}  A point $x \in  \AA$ lies in
$F_\AA$ if and only if 
$\gg_{x,r} = \gg_{F_\AA}$ and $\gg_{x,r^+} = \gg_{F_\AA}^+${\rm .}  
\endproclaim 
 
Thanks to a suggestion of Jiu-Kang~Yu, the proof below is far more
elegant than the original.

\demo{Proof}
The $r$-facet in $\AA$ to which $x$ belongs is
completely determined 
by the three sets
$$
\{\psi \in \Psi(\AA) \, | \, \psi(x) > r \} \qquad
\{\psi \in \Psi(\AA) \, | \, \psi(x) = r \} \qquad
\{\psi \in \Psi(\AA) \, | \, \psi(x) < r \}.
$$
These three sets are, in turn, completely determined by $\gg_{x,r}$
and $\gg_{x,r^+}$.
\enddemo

\specialnumber{3.1.5}\proclaim{Lemma} \label{lem:openone}
If $y \in \AA${\rm ,} then the union of all $r$\/{\rm -}\/facets of $\AA$ which contain $y$ in
their closure is an open neighborhood of $y$ in $\AA${\rm .}
\endproclaim

3.2. {\it Generalized $r$\/{\rm -}\/facets}.

\vglue8pt {\it Definition} {3.2.1}.
For $x \in \BB(G)$, define 
$$ F^*(x) := \{ y \in \BB(G) \,| \, \gg_{x,r} = \gg_{y,r} \hbox{ and }
\gg_{x,r^+} = \gg_{y,r^+} \}.$$

 {\it Definition} 3.2.2.
$${\cal  F}(r) := \{F^*(x) \, | \, x \in \BB(G)\}.$$

{\it Definition} 3.2.3.
An element of ${\cal  F}(r)$ is called a {\it generalized
$r$-facet}.

\specialnumber{3.2.4}\numbereddemo{{R}emark} \label{rem:grf}
 1. If $x \in \BB(G)$, then for all $y \in F^*(x)$ we have $F^*(x)
= F^*(y)$. 
\begin{itemize}

\item[2.] Suppose that $x,y \in \BB(G)$.  We write $x \sim y$ if and 
only if $F^*(x) = F^*(y)$.   Then   
$$\BB(G) = \coprod_{x \in \BB(G) / \sim} F^*(x) = \coprod_{F^* \in {\cal  F}(r)}
F^*.$$ 
\item[3.] For $x \in \BB(G)$ and $g \in G$ we have $gF^*(x) =
F^*(gx)$. 
\item[4.]  
If $F^* \in {\cal  F}(r)$ and $\AA$ is an apartment of
$\BB(G)$ such that $F_\AA = \AA \cap F^* \neq \emptyset$, then it follows from
Lemma~\ref{lem:rfacec} that $F_\AA$ 
is an $r$-facet of $\AA$. 
\item[5.] If $F^* \in {\cal  F}(r)$, then $F^*$ is a nonempty and convex
subset of $\BB(G)$. 
\end{itemize}
\enddemo

\specialnumber{3.2.5}\proclaim{Lemma} \label{lem:needed}
$F^* \in {\cal  F}(r)$ if and only if $F^* \in {\cal  F}(-r)${\rm .}
\endproclaim

{\it Proof}.
This follows from Remarks~3.1.2 and
3.2.4 (4). 
\hfill\qed

\specialnumber{3.2.6}\proclaim{Lemma} \label{lem:rfaceb}
If $x \in \BB(G)$ and $\AA$ is an apartment in $\BB(G)$ such that
$F_\AA = F^*(x) \cap \AA \neq \emptyset${\rm ,} then for all $y \in F_\AA$
we have 
$$F^*(x) = {G_y} F_{\AA}.$$
\endproclaim

\vglue-12pt
{\it Proof}.
Fix $y \in F_\AA$.

``$\subset$'': Suppose $z \in F^*(x)$.  Then there exists an $h \in G_y$
such that $hz \in \AA$.  Note that 
$$\gg_{hz, r} = \lsup{h} \gg_{z,r} = \lsup{h}\gg_{x,r} = \lsup{h}\gg_{y,r} =
\gg_{hy,r} = \gg_{y,r} = \gg_{x,r}$$
and similarly $\gg_{hz, r^+} = \gg_{x,r^+}$. 
Thus $hz \in \AA \cap F^*(x) = F_\AA$, and so $z \in {G_y} F_\AA$.

``$\supset$'':  Suppose $z \in F_\AA$ and $h \in G_y$.  We have
$$\gg_{hz, r} = \lsup{h} \gg_{z,r} = \lsup{h}\gg_{y,r} =
\gg_{hy,r} = \gg_{y,r} = \gg_{x,r}$$ 
and similarly 
$\gg_{hz, r^+} = \gg_{x,r^+}$. 
Thus $hz \in F^*(x)$.
\hfill\qed

\specialnumber{3.2.7}\proclaim{{C}orollary} \label{cor:rfacebdd}
If $F^* \in {\cal  F}(r)${\rm ,} then the image of $F^*$ in
$\BB^{\rm red}(G)${\rm ,} the reduced Bruhat\/{\rm -}\/Tits building{\rm ,} is
bounded{\rm .} 
\endproclaim

\specialnumber{3.2.8}\proclaim{Lemma} \label{lem:normstab}
For $x \in \BB(G)$ we have
$$N_G(\gg_{x,r}) \cap N_G( \gg_{x,r^+}) = \stab_G \bigl(F^*(x)\bigr).$$
\endproclaim

\vglue-24pt
{\it Proof}.
Let $F^* = F^*(x)$.

We have $\stab_G(F^*) \subset N_G(\gg_{x,r}) \cap N_G(\gg_{x,r^+})$.

Suppose $n \in N_G(\gg_{x,r}) \cap N_G(\gg_{x,r^+})$.  Fix $z \in
F^*$.  Let $\AA$ be an 
apartment of $\BB(G)$ containing $x$ and $nz$.  Let $F_\AA = \AA \cap
F^*$ ( $\neq \emptyset$). If $y \in \AA$ such that $\gg_{y,r} =
\gg_{x,r}$ and $\gg_{y,r^+} = \gg_{x,r^+}$, then from
Lemma~\ref{lem:rfacec} we have $y \in 
{F_\AA}$. 
Thus, since $\gg_{nz,r} = \lsup{n} \gg_{z,r} = \lsup{n} \gg_{x,r} =
\gg_{x,r}$ and similarly $\gg_{nz,r^+} = \gg_{x,r^+}$, we have $nz  \in
{F_\AA} \subset {F^*}$.  Since 
$z$ was arbitrary, we have $nF^* \subset {F^*}$.
\hfill\qed

\specialnumber{3.2.9}\proclaim{Lemma}  \label{lem:crface}
If $F^* \in {\cal  F}(r)$ and $\AA$ is an apartment in $\BB(G)$ such
that $F_\AA = F^* \cap \AA \neq \emptyset${\rm ,} then 
$$\overline{F_\AA} =
\overline{F^*} \cap \AA.$$ 
\endproclaim

\vglue-12pt
{\it Proof}.
Suppose $F^* \in {\cal  F}(r)$ and $\AA$ is an apartment in $\BB(G)$ such
that $F_\AA = F^* \cap \AA \neq \emptyset$.
It is enough to show that $\overline{F^*} \cap \AA \subset
\overline{F_\AA}$. 

Suppose $x \in \overline{F^*} \cap \AA$.  
Let $\{x_n\}$ be a sequence in $F^*$ which converges 
to $x$.   Fix $y \in F_\AA$.

By choosing a subsequence of $\{x_n\}$, we may assume
that for each $n \in \N$ there exists a zero-alcove $C_n$ such that
$x_n$ and $x$ both live in 
$\overline{C_n}$.  We may also assume that $\dist(x_n,x) < 1/n$ for
all $n \in \N$.  For the remainder of this
paragraph, fix $n \in \N$. 
Let $\AA_n$ be an apartment in
$\BB(G)$ containing both $C_n$ and $y$.  Since $x$ and $y$ both lie in
$\AA_n \cap \AA$, there exists
$g_n \in G$
such that $g_n$ fixes both $x$ and $y$ and $g_n \AA_n = \AA$.  Since
$g_n$ fixes $x$, we have
$$\dist(g_n x_n,x) = \dist(x_n,x) < 1/n.$$
Since $g_n$ fixes $y$, it follows from Lemma~\ref{lem:normstab} that
$g_n \in \stab_G(F^*)$.  Since $g_n \AA_n = \AA$ and $g_n \in
\stab_G(F^*)$, we have $g_n x_n \in F^* \cap \AA = F_\AA$.

Consequently, the sequence $\{g_n x_n\}$ in $F_\AA$
converges to $x$.  Thus $x \in \overline{F_\AA}$.
\phantom{great}\hfill\qed
\vglue8pt

 {\it Definition} 3.2.10. For $F^* \in {\cal  F}(r)$ and $\delta > 0$, define
$$F^*(\delta) := \{ x\in \overline{F^*} \, | \, \dist(x,z) \geq \delta
\, \hbox{ for all } \, z \in \overline{F^*} \smallsetminus F^* \}.$$

\specialnumber{3.2.11}\proclaim{Lemma} \label{lem:deltaok}
Suppose $F^* \in {\cal  F}(r)$ and $\delta > 0${\rm .} We have that
$F^*(\delta)$ is a convex{\rm ,} closed{\rm ,} and  
$\stab_G(F^*)$\/{\rm -}\/invariant subset of $\BB(G)${\rm .} 
Moreover{\rm ,} $F^*(\delta)$ is a nonempty subset of $F^*$ if and only if there exists an
apartment $\AA$ in $\BB(G)$ such that the subset of $F_\AA = F^* \cap
\AA$ defined by 
$$F_{\AA}(\delta) = \{ x \in F_\AA \, | \, \dist(x,z) \geq \delta
\, \hbox{ for all } \, z \in \overline{F_\AA} \smallsetminus F_\AA \}$$
is nonempty{\rm .}
\endproclaim

{\it Proof}.
$F^*(\delta)$ is a closed and $\stab_G(F^*)$-invariant subset of $\BB(G)$.
We now consider the last statement of the lemma.

For all apartments $\AA$ of $\BB(G)$ we have $F^*(\delta) \cap \AA
\subset F_\AA(\delta)$.  Thus, if $F^*(\delta)$ is nonempty, then
there exists an apartment $\AA$ in $\BB(G)$ such that $F_\AA(\delta)
\neq \emptyset$.

We will show that if there is an  apartment $\AA$ in $\BB(G)$
such that $F_\AA(\delta) 
\neq \emptyset$, then
\begin{equation} \label{equ:deltaversion}
G_y F_\AA(\delta) = F^*(\delta) 
\end{equation}
for all $y \in F_\AA(\delta)$. This implies that if $F_\AA(\delta)
\neq \emptyset$, then
$F^*(\delta) \neq \emptyset$.

Suppose $\AA$ is an apartment in $\BB(G)$ such that $\emptyset
\neq  F_\AA(\delta) \subset F_\AA = F^* \cap \AA$.  Fix $w \in
F_\AA(\delta)$.   

We first show that 
$G_w F_\AA(\delta) \subset F^*(\delta)$.  Since $G_w \leq
\stab_G(F^*)$, we have that $F^*(\delta)$ is $G_w$-invariant.  Thus, it
will be enough to show that $F_\AA(\delta) \subset F^*(\delta)$.  Fix
$x \in F_\AA(\delta)$.  Suppose $z \in \overline{F^*} \smallsetminus
F^*$.  Choose an apartment $\AA_z$ such that $x$ and $z$ both belong
to $\AA_z$.  There exists a $g \in G_x$ such that $g \AA_z = \AA$.
Since $g \in G_x$, it follows from Lemma~\ref{lem:normstab} that $g
\in \stab_G(F^*)$.  From Lemma~\ref{lem:crface} we have $gz \in
(\overline{F^*} \smallsetminus F^*) \cap \AA =
\overline{F_\AA} \smallsetminus F_\AA$.
Thus
$$\delta \leq \dist(x,gz) = \dist (x,z).$$
Since $z$ was arbitrary, we have $x \in F^*(\delta)$.

We now show that 
$G_w F_\AA(\delta) \supset F^*(\delta)$.  Fix $z \in F^*(\delta)$.
From Lemma~\ref{lem:rfaceb} there exists $k \in G_w$ such that $kz \in
F_\AA$.  Since $k \in G_w \subset \stab_G(F^*)$, we have $kz \in F_\AA
\cap kF^*(\delta) = F_\AA \cap F^*(\delta) \subset F_\AA(\delta).$  
Thus, equation~(\ref{equ:deltaversion}) is valid.

It remains to see that $F^*(\delta)$ is convex.  If $F^*(\delta)$ is empty,
there is nothing to prove.  So suppose $F^*(\delta)$ is nonempty.
Then there exists an apartment $\AA$ in $\BB(G)$ such that
$F_\AA(\delta)$ is 
nonempty.
Suppose $x,z
\in F^*(\delta)$.  Fix $w \in F_\AA(\delta)$.
From~(\ref{equ:deltaversion}) there exists $k \in 
G_w$ such that $kx \in F_\AA(\delta)$.  Since $k \in G_w \leq
\stab_G(F^*)$, we have $kz \in F^*(\delta)$.  Thus, another
application of~(\ref{equ:deltaversion}) shows that there exists $k_1
\in G_{kx}$ such that 
$k_1kz \in F_\AA(\delta)$.  As $F_\AA(\delta)$ is convex, we have 
$$[k_1 k x, k_1 k z] \subset F_\AA(\delta) \subset F^*(\delta).$$
Since $F^*(\delta)$ is
$\stab_G(F^*)$-invariant, we have $[x,z] \subset F^*(\delta)$. 
\hfill\qed\vglue8pt

 {\it Definition} 3.2.12. For $F^* \in {\cal  F}(r)$, define
$$C(F^*) :=\left\{ y \in F^* \, | \begin{array}{rl} & \hbox{for all apartments $\AA$ of
$\BB(G)$ for which}\\
&\hbox{ $\AA \cap F^* \neq \emptyset$ we have $y \in \AA$}.\end{array}\right\} $$

\specialnumber{3.2.13}\proclaim{{C}orollary} \label{cor:notem}
If $F^* \in {\cal  F}(r)$, then $C(F^*) \neq \emptyset${\rm .} 
\endproclaim

{\it Proof}.
Without loss of generality, we suppose that $\bGn$ is semisimple.
Let $N = \stab_G(F^*)$.  From Corollary~\ref{cor:rfacebdd} we
have that $F^*$ is bounded in $\BB(G)
=\BB^{\rm red}(G)$.  Thus it follows from
~\cite[\S2.2.1]{corvallis} that $N$ is  
a bounded subgroup of~$G$.

If $F^*$ consists of a single point, there is nothing to prove.  So we
suppose that $F^*$ is not a point.  Let $\AA$ be an
apartment in $\BB(G)$ such that $F_\AA = \AA \cap F^* \neq \emptyset$.
It follows from Lemma~\ref{lem:rfaceb} that $\dim F_\AA > 0$.  Thus,
there exists $\delta > 0$ such that the set
$$\{ x \in F_\AA \, | \, \dist(x,z) \geq \delta
\, \hbox{ for all } \, z \in \overline{F_\AA} \smallsetminus F_\AA \}$$
is nonempty.  From
Lemma~\ref{lem:deltaok} there exists $\delta > 0$ such that
$F^*(\delta)$ is a nonempty, convex, closed, $N$-stable subset of
$\BB(G)$.  Consequently, there exists a $y \in F^*(\delta) \subset F^*$
such that $ny =y$ for all $n \in
N$~\cite[{Proposition}~3.2.4]{bruhat-tits:one}.  

We now show that $y \in C(F^*)$.  Suppose $\AA'$ is an apartment of
$\BB (G)$ such that 
$F_{\AA'} = \AA' \cap 
F^* \neq \emptyset$.  Choose $z \in F_{\AA'}$.  From
Lemma~\ref{lem:rfaceb} we have $G_z F_{\AA'} = F^*$.  However, $G_z
\subset N$ from 
Lemma~\ref{lem:normstab}.  Thus $G_zy = y$.  Consequently, we must
have $y \in F_{\AA'} \subset \AA'$. 
\hfill\qed

\specialnumber{3.2.14}\proclaim{{C}orollary} \label{cor:rfacedim}
If $\AA_i$ ($i = 1,2$) are two apartments of $\BB(G)$ and $F^* \in
{\cal  F}(r)$ 
such that $F_{\AA_i} = F^* \cap \AA_i \neq \emptyset${\rm ,} then 
$\dim A(F_{\AA_1}, \AA_1) = \dim A(F_{\AA_2},\AA_2)${\rm .}
\endproclaim

\specialnumber{3.2.15}\proclaim{Lemma} \label{lem:clwork}
If $F_i^* \in {\cal  F}(r)$ ($i = 1,2$) such that $F_1^* \cap
\overline{F_2^*} \neq \emptyset${\rm ,} then $F_1^* \subset \overline{F_2^*}${\rm .} 
\endproclaim

{\it Proof}.
Fix $y_i \in C(F_i^*)$.  We first show that $y_1 \in F^*_1 \cap
\overline{F_2^*}$.   Let $z \in F_1^* \cap \overline{F_2^*}$.  Choose an
apartment $\AA$ containing $z$ and $y_2$.  Let $F_{i,\AA} = F^*_i
\cap \AA$.  Since\break $F_{1,\AA} \neq \emptyset$, we have $y_1 \in
F_{1,\AA}$.  We also have 
$z \in \overline{F_2^*} \cap \AA = 
\overline{F_{2,\AA}}$
from Lemma~\ref{lem:crface}.
Now $z \in F_{1,\AA} \cap \overline{F_{2,\AA}}$ so $F_{1,\AA} \subset
\overline{F_{2,\AA}}$ since these are both $r$-facets of $\AA$. Thus
$y_1 \in  \overline{F_2^*}$.

Suppose $w \in F_1^*$.  Let $\AA'$ be an apartment containing $w$ and
$y_2$.  Let $F_{i,\AA'} = F^*_i \cap \AA'$.  From the previous
paragraph we have $y_1 \in
F_{1,\AA'}$ and $y_1 \in 
\overline{F^*_2} \cap \AA' = \overline{F_{2,\AA'}}$.  Since
$F_{2,\AA'}$ and $F_{1,\AA'}$ are both $r$-facets 
of $\AA'$, we have $w \in F_{1,\AA'} \subset \overline{F_{2,\AA'}}$. 
\hfill\qed
\vglue4pt

  Thanks to
Corollary~\ref{cor:rfacedim} the following definition makes sense.
\vglue4pt
{\it Definition} 3.2.16. Suppose $F^* \in {\cal  F}(r)$.  Let $\AA$ be an apartment in
$\BB(G)$ such that $\AA \cap F^* \neq \emptyset$.  We define
$$\dim F^* := \dim A(F^* \cap \AA, \AA).$$
\vglue4pt

Moreover, it follows from Lemma~\ref{lem:clwork} that $\overline{F^*}$ is
the disjoint union of $F^*$ and generalized $r$-facets which meet
$\overline{F^*}$ and have
dimension strictly smaller than that of $F^*$.

\specialnumber{3.2.17}\proclaim{Lemma}
If $F_i^* \in {\cal  F}(r)$ ($i = 1,2$) such that $F_1^* \neq F_2^*$
and $F_1^* \subset
\overline {F_2^*}${\rm ,} then for fixed $y_i \in C(F_i^*)$ there exists an $x_2
\in F_2^*$ such that 
 \vglue4pt
{\rm 1.} $G_{x_2} \subset G_{y_1}$ and
 \vglue4pt
{\rm 2.} $x_2 \in (y_1, y_2]${\rm .}
\endproclaim
 
\demo{Proof}
Choose an apartment $\AA$ containing $y_1$ and $y_2$.  Let $F_{i,\AA}
= \AA \cap F_i^*$.  We have
$F_{1,\AA} = F_1^* \cap \AA \subset \overline{F_2^*} \cap \AA =
\overline{F_{2,\AA}}$.
Let $${\cal  F}(y_1,0) = \{ H \in {\cal  F}(0) \, | \, H \subset
\AA \hbox{ and } y_1 \in \overline{H} \}.$$  
Note that $\bigcup_{H \in {\cal  F}(y_1,0)} H$ is an open
neighborhood of $y_1$ in $\AA$.  Consequently, there exists an $H \in
{\cal  F}(y_1, 0)$ such that $H \cap (y_1,y_2] \neq \emptyset$. 

Choose $x_2 \in H \cap (y_1, y_2] \subset F_{2,\AA}$.  We have
$G_{x_2} = G_H \subset G_{y_1}$. 
\enddemo

 {\it Definition} 3.2.18.
Suppose $F^* \in {\cal  F}(r)$.  Fix $x \in F^*$.  We define
$$\gg_{F^*} := \gg_{x,r}$$
and
$$\gg_{F^*}^+ := \gg_{x,r^+}.$$
\vglue6pt

Sometimes, to avoid confusion, we denote $\gg_{F^*}$ by $\gg_{F^*,r}$
and $\gg_{F^*}^+$ by $\gg_{F^*,r^+}$.  We now present a
corollary to the proof of Lemma~\ref{lem:clwork}. 

\specialnumber{3.2.19}\proclaim{{C}orollary} \label{cor:closure} 
Suppose $F_i^* \in {\cal  F}(r)$ for $i = 1,2${\rm .}  If $F_1^* \subset
\overline{F_2^*}${\rm ,} then  
$$\gg_{F_1^*}^+ \subset \gg_{F_2^*}^+ \subset \gg_{F_2^*} \subset
\gg_{F_1^*}.$$ 
\endproclaim

{\it Proof}.
Choose $y_i \in C(F_i^*)$.  Let $\AA$ be an apartment in $\BB(G)$
containing $y_1$ and $y_2$.  Let $F_{i,\AA} = F_i^* \cap \AA$.  From
the proof of Lemma~\ref{lem:clwork}, we have $F_{1,\AA} \subset
\overline{F_{2,\AA}}.$  We then have $\gg_{F_{1,\AA}}^+ \subset
\gg_{F_{2,\AA}}^+ \subset \gg_{F_{2,\AA}} \subset \gg_{F_{1,\AA}}.$ 
\hfill\qed

\specialnumber{3.2.20}\proclaim{Lemma}  \label{lem:opengen}
If $y \in \BB(G)${\rm ,} then
the union of all generalized $r$\/{\rm -}\/facets that contain $y$ in their
closure is an open neighborhood of $y$ in  $\BB(G)${\rm .}
\endproclaim  

\demo{Proof}
Fix $y \in \BB(G)$ and an apartment $\AA$ in $\BB(G)$ which contains $y$.  Let
$H_\AA$ denote the union of all $r$-facets of $\AA$ which contain $y$
in their closure.  From Lemma~\ref{lem:openone} the set $H_\AA$ is
an open neighborhood of $y$ in $\AA$.  Fix $\varepsilon > 0$ so that if
$x \in \AA$ and $\dist(x,y) < \varepsilon$, then $x \in H_\AA$.

Let $H$ denote the union of  all generalized $r$-facets that contain
$y$ in their closure.  We will show that the ball in $\BB(G)$ of radius
$\varepsilon$  centered around $y$ is contained in $H$.
Fix $z \in \BB(G)$ such that 
$\dist(z,y) < \varepsilon$.
There exists $g \in G_y$ such that $gz \in \AA$.  Since $\dist(gz,y) =
\dist (gz, gy) = \dist(z,y) < \varepsilon$, we have $gz \in H_\AA$.
Since $gz \in H_\AA$, there
exists $F^* \in {\cal  F}(r)$ such that $y \in \overline{F^* \cap
\AA}$ and $gz \in F^* \cap
\AA$.  Thus we have $gz \in {F^*}$ and $y \in \overline{F^*}$.  Since
$y = g\inv y \in g\inv \overline{F^*} = \overline{g\inv F^*}$, we
conclude that $z \in g\inv F^* \subset H$.
\enddemo

3.3. {\it Associativity}.
In this subsection we introduce an equivalence relation on the
elements of ${\cal  F}(r)$ which is a generalization of  the concept of
``associate'' found in~\cite{moy-prasad:jacquet}, \cite{moy-prasad:K-types}
(see Remark~\ref{rem:associate}).

\specialnumber{3.3.1}\numbereddemo{Definition}Suppose $F^* \in {\cal  F}(r)$ and $\AA$ is an apartment in
$\BB(G)$.  We define
$A(\AA,F^*) := A(F^* \cap \AA, \AA)$.
\enddemo

 {\it Definition} 3.3.2. Two generalized $r$-facets $F^*_1$ and $F^*_2$ are {\it strongly
$r$-associated} if for all apartments $\AA$ such that  $F_1^* \cap \AA
\neq \emptyset$ and $F_2^* \cap \AA
\neq \emptyset$, we have
$$A(\AA, F_1^*) = A(\AA,F_2^*).$$

\specialnumber{3.3.3}\proclaim{Lemma}
Two generalized $r$\/{\rm -}\/facets $F_1^*, F^*_2 \in {\cal  F}(r)$ are
strongly\break $r$\/{\rm -}\/associated if and only if there exists an apartment $\AA$
such that 
$\emptyset \neq A(\AA, F_1^*) = 
A(\AA,F_2^*)$\/{\rm .}\/
\endproclaim

\demo{Proof}
``$\Rightarrow$'':  This follows from the definition.

``$\Leftarrow$'':  Choose $x_i \in C(F_i^*)$ for $i = 1,2$.  Recall that
for an apartment $\AA'$ of $\BB(G)$ we have $\AA' \cap F_i^* \neq
\emptyset$ if and only if $x_i \in \AA'$.  Suppose
$\AA' \cap F_1^* \neq \emptyset$ and $\AA' \cap F_2^* \neq \emptyset$.  There
exists a  $g \in G$ such that $g$ fixes $\AA \cap \AA'$
point-wise and $g\AA = \AA'$.  Thus $gx_1 = x_1$  and
$gx_2 = x_2$.  This implies that  $g
\in \stab_G(F_1^*) \cap \stab_G(F_2^*)$ and
\begin{eqnarray*}
A(\AA',F_1^*) &=& A(g\AA,gF_1^*) = g A(\AA,F_1^*)\\
&=& gA (\AA, F_2^*) = A(g\AA, gF_2^*) \\
&=&A(\AA', F_2^*). \\
\noalign{\vskip-36pt}
\end{eqnarray*}
\enddemo

\specialnumber{3.3.4}\numbereddemo{Definition}Two generalized $r$-facets $F_1^*$ and $F_2^*$ are
{\it $r$-associated} if there exists a $g\in G$ such that $F_1^*$ and
$gF_2^*$ are strong $r$-associates.
\enddemo

\specialnumber{3.3.5}\numbereddemo{{R}emark} \label{rem:associate}
If $F_1^*, F_2^* \in {\cal  F}(0)$ are $0$-associated, then the 
parahoric subgroups $G_{F_1^*,0}$ and $G_{F_2^*,0}$ are associate in
the sense of~\cite{moy-prasad:K-types}.
\enddemo

\demo{Example {\rm 3.3.6}}
In Figure~1, we have represented a $0$-alcove in the
building of ${\bf SL}_3(k)$ (resp., ${\bf G}_2(k)$).  The edges
identified with 
hatch marks are $0$-associates; none of the remaining pictured $0$-facets are
$0$-associated.
\enddemo

 \figin{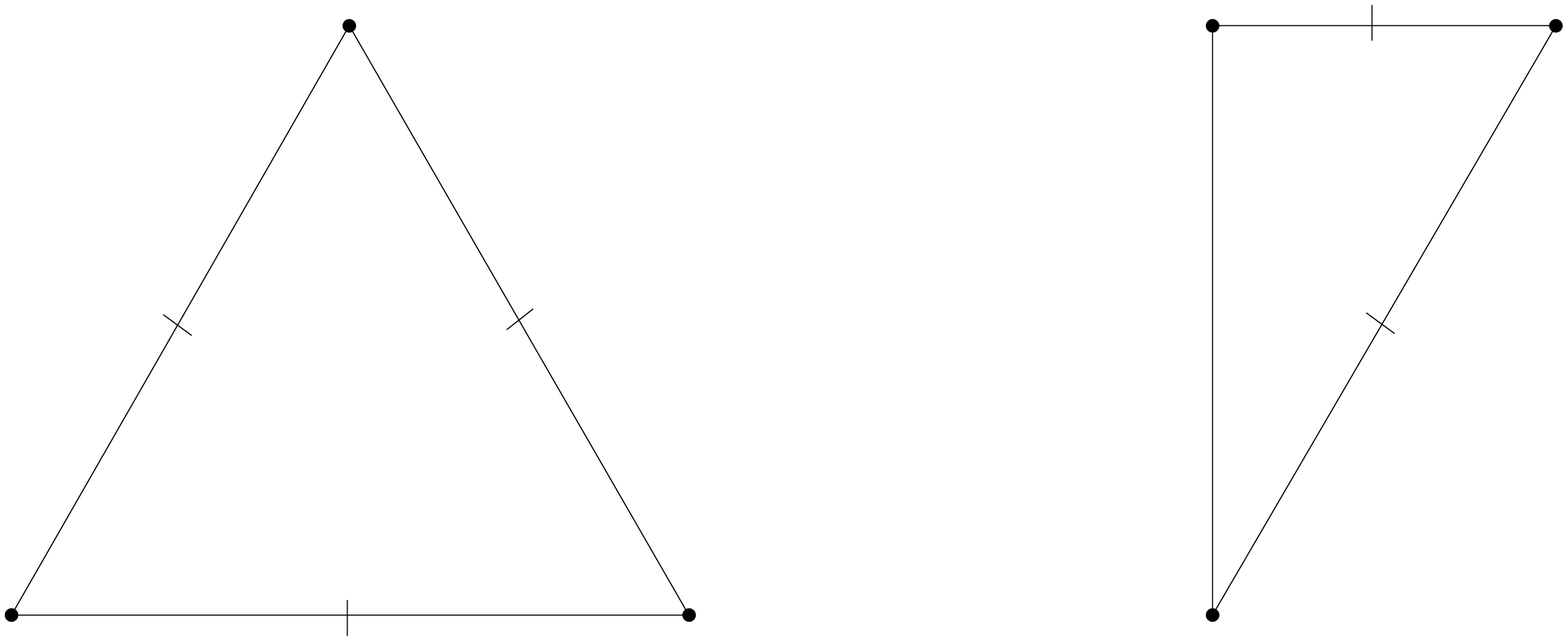}{300}
 \begin{quote}
 Figure 1. {Associates in $0$-alcoves for ${\bf SL}_3(k)$ (resp.,
 ${\bf G}_2(k)$).} 
 \end{quote}

\specialnumber{3.3.7}\proclaim{Lemma} \label{lem:rassocequiv}
 $r$-associativity is an
equivalence relation on ${\cal  F}(r)$. 
\endproclaim

{\it Proof}.
For two generalized $r$-facets $F_1^*$ and $F_2^*$, we
write $F_1^* \sim F_2^*$ if and only if $F_1^*$ and $F_2^*$ are
$r$-associated.
The relation is reflexive and symmetric.  We now show that it
is transitive.  

Suppose $F_1^*, F_2^*, F_3^* \in {\cal  F}(r)$ 
such that $F_1^* \sim F_2^*$ and $F_2^* \sim F_3^*$.
There exist $g_2, g_3 \in G$ and apartments $\AA_{12}, \AA_{23}$ in
$\BB(G)$ such that
$$A(\AA_{12},F_1^*) = A(\AA_{12}, g_2 F_2^*) \neq \emptyset$$
and
$$A(\AA_{23},F_2^*) = A(\AA_{23}, g_3 F_3^*) \neq \emptyset.$$
Let $z \in C(F_2^*)$.  Then $z \in g_2\inv \AA_{12} \cap \AA_{23}$ and
so there exists $h \in G_z \subset \stab_G (F_2^*)$ such that
$hg_2\inv \AA_{12} = \AA_{23}$.  We have
\begin{eqnarray*}
A(\AA_{12},F_1^*) &=& A(\AA_{12},g_2 F_2^*) = g_2 A(g_2\inv \AA_{12},
F_2^*)\\
&= &g_2 h\inv A(\AA_{23},
F_2^*) = g_2 h\inv A( \AA_{23}, g_3 F_3^*) \\
&=& A(\AA_{12}, g_2h\inv g_3 F_3^*).
\\
\noalign{\vskip-36pt}
\end{eqnarray*}
\hfill\qed

\specialnumber{3.3.8}\numbereddemo{{R}emark}
${\cal  F}(r)/ \! \sim$ is finite.
\enddemo

3.4. {\it Some finite\/{\rm -}\/dimensional vector spaces}.

\specialnumber{3.4.1}\numbereddemo{Definition}
For $x \in \BB(G)$ denote the finite-dimensional
$\ff$-vector space $\gg_{x,r} / \gg_{x,r^+}$ by $V_{x,r}$.
\enddemo

{\it Definition} 3.4.2.
If $F^* \in {\cal  F}(r)$ and $x \in F^*$, then $V_{F^*} := V_{x,r}$.
 
\specialnumber{3.4.3}\numbereddemo{Definition}
If $\AA$ is an apartment in $\BB(G)$,  $F_\AA$ is an $r$-facet of
$\AA$, and $x \in F_\AA$, then $V_{F_\AA} := V_{x,r}$.
\enddemo

3.5. {\it A natural identification}.
In this subsection we show that if $F_1^*, F_2^* \in {\cal  F}(r)$ are
strongly $r$-associated, then we can naturally identify $V_{F_1^*}$ with
$V_{F_2^*}$.  Moreover, we show that these two spaces have the same
orbit structure under this identification.

\specialnumber{3.5.1}\proclaim{Lemma} \label{lem:natident}
If $F_1^*, F_2^* \in {\cal  F}(r)$ are strongly $r$\/{\rm -}\/associated{\rm ,} then
the natural map
$$\gg_{F_1^*} \cap \gg_{F_2^*} \rightarrow V_{F_i^*}$$
is surjective with kernel $\gg_{F_1^*}^+ \cap \gg_{F_2^*} = \gg_{F_1^*}
\cap \gg_{F_2^*}^+ = \gg_{F_1^*}^+ \cap \gg_{F_2^*}^+${\rm .}
\endproclaim

{\it Proof}.
Choose an apartment $\AA$ in $\BB(G)$ for which $F_{i,\AA} = F_i^*
\cap \AA \neq \emptyset$ for $i = 1,2$.  If $\psi \in \Psi(\AA)$ such
that $\psi|_{F_{i,\AA}} = r$, then $A(\AA,F_i^*) \subset H_{\psi-r}
\subset \AA$.  Thus, since $F_{1,\AA}$ and $F_{2,\AA}$ are open in
$A(\AA,F_2^*) = 
A(\AA,F_1^*)$, we have

\centerline{${\displaystyle\psi|_{F_{1,\AA}} = r \hbox{ if and only if } \psi|_{F_{2,\AA}} = r}$}

\noindent 
 for all $\psi \in \Psi(\AA)$.
The lemma follows.
\hfill\qed

\specialnumber{3.5.2}\numbereddemo{{R}emark}
From Lemma~\ref{lem:natident}, we obtain a bijective identification of
$V_{F_1^*}$ with $V_{F_2^*}$.  We write
$$V_{F_1^*} \ident V_{F_2^*}$$
for this identification.  More generally, we will use the ``$\ident$''
notation whenever two objects are to be identified via this natural bijection.
\enddemo

{\it Definition} 3.5.3. If $F^* \in {\cal  F}(r)$ and $x \in F^*$, then the image of $G_x$
in $\aut_{\ff}(V_{F^*})$ is denoted by $N_x(F^*)$.

\specialnumber{3.5.4}\proclaim{Lemma} \label{lem:identaut}
Suppose $F_i^* \in {\cal  F}(r)$ and $x_i \in F_i^*$ for $i = 1,2${\rm .} If
$F_1^*$ and $F_2^*$ are strongly $r$\/{\rm -}\/associated{\rm ,} then $N_{x_i}(F_i^*)$ is
the image of $G_{x_1} \cap G_{x_2}$ in $\aut_{\ff}(V_{F_i^*})$ for $i=1,2${\rm .}
Moreover 
$$N_{x_1} (F_1^*) \ident N_{x_2}(F_2^*)$$
under the identification of $V_{F_1^*}$ with $V_{F_2^*}$ introduced
above{\rm .} 
\endproclaim

{\it Proof}.
We have that $V_{F_i^*}$ is the image of $\gg_{x_1,r} \cap
\gg_{x_2,r}$ in $V_{F_i^*}$(= $\gg_{x_i,r}/\gg_{x_i,r^+}$) for $i=1,2$.
Let $\AA$ be an apartment in $\BB(G)$ 
containing $x_1$ and $x_2$.  Suppose $\psi \in \Psi(\AA)$ such that
$\psi(x_1) = 0$ and the image of $U_{\psi}$ in $\aut_{\ff}(V_{F_1^*})$ is
nontrivial.  Since the image of $U_\psi$ is nontrivial, there exist
$X \in \gg_{x_1,r} \cap \gg_{x_2,r}$ and $g \in U_\psi$ such that
$\lsup{g}X \neq X \hbox{ mod } (\gg_{x_1,r^+} \cap \gg_{x_2, r^+})$.

We now show that $\psi(x_2) = 0$.

If $\psi(x_2) > 0$,  then $\lsup{g}X = X \hbox{ mod }
\gg_{x_2, r^+}$.  Since $\lsup{g}X - X \in \gg_{x_1,r}$, from
Lemma~\ref{lem:natident} we have $\lsup{g}X = X \hbox{ mod } (\gg_{x_1,r^+}
\cap \gg_{x_2, r^+})$.  We therefore conclude that $\psi(x_2) \leq 0$.

If $\psi(x_2) < 0$, then we let $v$ denote the vector $(x_2
- x_1)$.  For all $\varepsilon \in \R$ we have $x_1 + \varepsilon \cdot
v \in
\AA$.  Consider the function 
$$f_v \colon \R \rightarrow \R$$
which sends $\varepsilon$ to $\psi(x_1 + \varepsilon \cdot v)$.  Note that
$f_v(0) = 0$ and $f_v(1) < 0$.  Since $\psi$ is affine, we have
that $f_v(\varepsilon ) > 0$ for all $\varepsilon < 0$.  Since $F_1^*
\cap \AA$ is open in $A(\AA , F_1^*) = A( \AA, F_2^*)$ and $x_1 + \R
\cdot v$ is an affine subspace of $A(\AA,F_1^*)$, there exists an
$\varepsilon < 0$ such that $x_1+ \varepsilon \cdot v \in F_1^* \cap
\AA$.  Thus for some $\varepsilon < 0$ we have  $x_1 + \varepsilon \cdot
v \in F_1^*$ and $U_\psi \subset
G_{x_1+ \varepsilon \cdot v}^+$.  Consequently, $g \in U_\psi$
acts trivially on $\gg_{F_1^*} \hbox{ mod } \gg^+_{F_1^*}$.

We therefore conclude that $\psi(x_2) = 0$.
Thus, if $g \in G_{x_1}$ has nontrivial image in
$\aut_{\ff}(V_{F_1^*})$, then it follows that we may
assume that $g \in G_{x_1} \cap G_{x_2}$.  The remainder of the lemma
now follows.
\hfill\qed\vglue8pt

From Lemma~\ref{lem:identaut} we have that
if $F^* \in {\cal  F}(r)$ and $x,y \in F^*$, then $N_x(F^*) =
N_y(F^*)$.  Therefore, the following definition makes sense.

\specialnumber{3.5.5}\numbereddemo{Definition}
If $F^* \in {\cal  F}(r)$ and $x \in F^*$, then define $N(F^*)
\subset \aut_{\ff}(V_{F^*})$ by
$$N(F^*) := N_x(F^*).$$
\enddemo

We can now restate Lemma~\ref{lem:identaut}.

\specialnumber{3.5.6}\proclaim{{C}orollary} \label{cor:paraaction}
If $F_1^*, F_2^* \in {\cal  F}(r)$ are strongly $r$\/{\rm -}\/associated{\rm ,} then
$N(F_1^*) \ident N(F_2^*)$. 
\endproclaim

3.6. {\it An equivalence relation on depth $r$ cosets}. 
In this subsection we introduce the set $I_r$ and an equivalence relation
on $I_r$.  The set $I_r$ parametrizes the set of all cosets of the
form $X+ \gg_{x,r^+}$ where $x \in \BB(G)$ and $X \in \gg_{x,r}$. 

\specialnumber{3.6.1}\numbereddemo{Definition}
$$I_r := \{(F^*,v) \, | \, F^* \in {\cal  F}(r) \hbox{ and } v \in
V_{F^*}\}.$$ 
\enddemo

We now introduce a relation on $I_r$.  Roughly speaking, two elements
$(F_1^*,v_1)$ and $(F_2^*,v_2)$ of $I_r$ are identified if
(1) $F_1^*$ and $F_2^*$ are $r$-associated, and (2) $v_1$ can then be
identified with a twist of $v_2$
(under the natural identification of the previous subsection).

\specialnumber{3.6.2}\numbereddemo{Definition}\label{def:irequiv}
For $(F_1^*,v_1)$ and $(F_2^*, v_2)$ in $I_r$ we write $(F_1^*,v_1)
\sim (F_2^*,v_2)$ if and only if there exist a $g \in G$ and an apartment $\AA$ in $\BB(G)$ such that
 
1. 
 $\emptyset \neq A(\AA, F_1^*) = A(\AA,gF_2^*)$ and

2. 
$\lsup{g}v_2 \ident v_1 \hbox{ \, in \, } V_{gF_2^*} \ident V_{F_1^*}$.
\enddemo

Here $\lsup{g}v_2$ has the obvious interpretation: if $X_2 \in
\gg_{F^*_2}$ is any lift of $v_2$,  then $\lsup{g}v_2$ denotes the
image of $\lsup{g}X_2$ in $V_{gF_2^*}$.

\specialnumber{3.6.3}\proclaim{Lemma}
The relation defined in Definition~{\rm \ref{def:irequiv}} is an equivalence relation on $I_r${\rm .}
\endproclaim

{\it Proof}.
The relation is reflexive.  We will show the relation is transitive;
once we do this, one can prove that the relation is symmetric in a
similar fashion.

We now show that the relation is transitive.  Suppose that
$(F_1^*,v_1), (F_2^*,v_2),\break (F_3^*,v_3) \in I_r$ such that $(F_1^*,v_1)
\sim (F_2^*,v_2)$ and $(F_2^*,v_2)
\sim (F_3^*,v_3)$.
Then there exist $g_2, g_3 \in G$ and apartments $\AA_{12}, \AA_{23}$
of $\BB(G)$ such that
$$\emptyset \neq A(\AA_{12}, F_1^*) = A(\AA_{12},g_2 F_2^*)$$
$$\emptyset \neq A(\AA_{23}, F_2^*) = A(\AA_{23},g_3 F_3^*),$$
 and
$$\lsup{g_2}v_2 \ident v_1 \hbox{ \, in \, } V_{g_2 F_2^*} \ident  V_{F_1^*}$$
$$\lsup{g_3}v_3 \ident v_2 \hbox{ \, in \, } V_{g_3 F_3^*} \ident  V_{F_2^*}.$$

We now wish to show that $(F_1^*,v_1) \sim (F_3^*,v_3)$.  We claim
that in
Definition~\ref{def:irequiv} the role of the pair ($g$,$\AA$) will be
played by ($h''h\inv g_2 g_3$, $\AA_{12}$) where $h \in G_{g_2 x_2}$ and
$h'' \in G_{h\inv g_2 g_3 x_3} \cap G_{h\inv g_2 x_2}$  will be
specified below.  
\vglue2pt
Fix $x_i \in C(F_i^*)$.  There exists an element $h \in G_{g_2 x_2}$
such that $h \AA_{12} = g_2 \AA_{23}$.  As in the proof of
Lemma~\ref{lem:rassocequiv} we have 
\begin{eqnarray*}
\emptyset \neq A(\AA_{12},F_1^*) &= &A(\AA_{12},g_2 F_2^*) = h\inv  A(h
\AA_{12}, h g_2 F_2^*)\\[6pt]
&=& h\inv g_2 A(  \AA_{23}, F_2^*) = h \inv g_2  A(\AA_{23},
g_3 F_3^*)\\[6pt]
&=& A(\AA_{12}, h\inv g_2 g_3 F_3^*).
\end{eqnarray*}
Arguing as in the proof of Lemma~\ref{lem:natident} we have that
$\gg_{F_1^*} \cap 
\gg_{g_2F_2^*} \cap \gg_{h\inv g_2 
g_3 F_3^*}$ surjects, under the natural map, onto $V_{F_1^*}$
(resp., $V_{g_2F_2^*}$, resp., $V_{h\inv g_2 g_3 F_3^*}$).  Choose $X
\in \gg_{F_1^*} \cap \gg_{g_2F_2^*} \cap \gg_{h\inv g_2
g_3 F_3^*}$ such that the image of $X$ in $V_{F_1^*}$ is $v_1$.

We have that the image of $X$ in $V_{g_2F_2^*}$ is $\lsup{g_2}v_2$.
Thus the image of $\lsup{g_2\inv}X$ in $V_{F_2^*}$ is $v_2$.  Since
$g_2\inv h g_2 \in G_{x_2}$, this implies that the image of
$\lsup{g_2\inv h}X = \lsup{(g_2\inv h g_2) g_2\inv}X$ in $V_{F_2^*}$
is $\lsup{g_2\inv h g_2}v_2$.   Note that $\lsup{g_2\inv h}X \in
\gg_{F_2^*} \cap \gg_{g_3 F_3^*}$.  Recall from
Corollary~\ref{cor:paraaction} that $N(F_2^*) \ident N(g_3F_3^*)$.
Thus, from Lemma~\ref{lem:identaut} there exists an $h' \in G_{g_3x_3}
\cap G_{x_2}$ such that 
$$\lsup{g_2\inv h g_2}v_2 \ident \lsup{h'g_3}v_3 \hbox{ \, in \, } V_{F_2^*}
\ident V_{g_3F_3^*}.$$
Thus, the image of $X$ in $V_{h\inv g_2g_3F_3^*}$ is
$$\lsup{h\inv g_2 h' g_3}v_3 = \lsup{h\inv (g_2 h' g_2\inv)g_2 g_3}v_3
= \lsup{h''h \inv g_2 g_3} v_3$$
where $h'' \in h\inv g_2 (G_{g_3x_3} \cap G_{x_2})g_2\inv h \subset
G_{h \inv g_2 g_3 x_3}$.  We have shown that
$$\emptyset \neq A(\AA_{12},F_1^*) = A(\AA_{12}, h \inv g_2 g_3 F_3^*) = A(\AA_{12},
h''h\inv g_2 g_3 F_3^*)$$
and
$$v_1 \ident \lsup{h'' h\inv g_2 g_3} v_3 \hbox{ \, in \, } V_{F_1^*} \ident
V_{h''h\inv g_2 g_3 F_3^*}.$$
So the relation is transitive.
\hfill\qed\vglue12pt

{\it {R}emark} 3.6.4.
If $\ff$ is finite, then $I_r/ \! \sim$ is finite.
 \pagebreak

\section{Jacobson-Morosov triples over $\ff$ and $k$} \label{sec:jacobson}

Fix $r \in \R$.  Much of the material in this section may be thought
of as a generalization of the material
in~\cite[\S3]{barbasch-moy:local}.  

In Section~4.3, we start with an $x \in \BB(G)$ and an
${\rmsl}_2(\ff)$-triple in $V_{x,-r}  
\times V_{x,0}  \times V_{x,r}$. From this data we  manufacture an
${\rmsl}_2(k)$-triple in $\gg$ which descends to our ${\rmsl}_2(\ff)$-triple. 

In Section~4.5 we perform this process in reverse.
That is, we start with an 
${\rmsl}_2(k)$-triple in $\gg$ and produce an 
$x \in \BB(G)$ such that our given ${\rmsl}_2(k)$-triple descends to an
${\rmsl}_2(\ff)$-triple in $V_{x,-r} \times V_{x,0}  \times V_{x,r}$.

\vglue4pt 4.1. {\it Degenerate cosets}.

\vglue4pt  {\it Definition} 4.1.1.
Suppose $F^* \in {\cal  F}(r)$.  An element $e \in V_{F^*}$ is
{\it degenerate} if and only if there exists a lift $E \in \gg_{F^*}$
of $e$ such that $E \in \NN$.

\proclaimtitle{Moy and Prasad}
\specialnumber{4.1.2}\proclaim{Lemma} \label{lem:degbenil}
Fix $F^* \in {\cal  F}(r)${\rm .} An element $e \in V_{F^*}$ is degenerate
if and only if zero is in the Zariski closure of \, $\lsup{\bfG_x}e$ for
all $x \in F^*${\rm .}
\endproclaim

{\it Proof}.
``$\Rightarrow$'': Fix $x \in F^*$.  Suppose $E \in \gg_{x,r} \cap \NN$
is a lift of $e$.  The desired conclusion follows
from~\cite[Proposition~4.3]{moy-prasad:K-types}.

``$\Leftarrow$'': This may also be derived
from~\cite{moy-prasad:K-types}.  We offer a slightly different proof.

We need to produce an $E \in \NN \cap \gg_{F^*}$ such that $E$ is a lift
of $e$.

Fix $x \in F^*$.  Let $\bS$ be a maximal $k$-split torus of $\bG$ such that $x
\in \AA\bigl(\bS(k)\bigr)$.  From~\cite{kempf:instability} 
there exists a 
one-parameter subgroup $\bar{\nu} \in \X_*^\ff(\bfG_x)$ such that $\lim_{t
\rightarrow 0} \lsup{\bar{\nu}(t)} e = 0$.  Let $\bfS$ be the maximal
$\ff$-split torus of $\bfG_x$ corresponding to $\bS$.  Since
maximal $\ff$-split tori are $\bfG_x(\ff)$-conjugate, there exist
$\bar{\mu} \in \X_*(\bfS)$ and $\bar{g} \in \bfG_x(\ff)$ such that
$\lim_{t \rightarrow 0} \lsup{\bar{\mu}(t) \bar{g}}e = 0$.  Let $\mu \in
\X_* (\bS)$ be the lift of $\bar{\mu}$ and let $g \in G_x$ be a lift of
$\bar{g}$.  Let $E' \in \gg_{x,r} = \gg_{F^*}$ be any lift of $e$.
We have $\lsup{g}(E' + \gg_{x,r^+}) = \lsup{g}(E' + \gg_{F^*}^+)
\subset \gg_{x+\varepsilon \cdot \mu, r^+}$ for all $\varepsilon$
sufficiently small and positive.  Consequently,
from~\cite{adler-debacker:kirilovn} we have $(E' +
\gg_{F^*}^+) \cap \NN 
\neq \emptyset$.  Choose $E$ in this intersection.
\hfill\qed\vglue9pt

4.2. {\it Some hypotheses}.  
  The statements below list properties
  which I require; no attempt has
  been made to produce a minimal list of hypotheses.
  If we 
  assume that $p$ is larger than some
  constant which can be determined by examining the absolute root
  datum of $\bGn$, then all of the hypotheses are valid.   In
  particular, if $\ff$ has characteristic zero, then the following
  hypotheses always hold.
  Where appropriate, I have identified references where a
  discussion about the conditions under which the hypothesis is
  valid may be found.\pagebreak

We begin by defining a finite-dimensional $\ff$-Lie algebra
$\overline{\gg}_x$.  Since we have fixed a uniformizer $\varpi$ for
$k$, for $s \in \R$ and $j \in \Z$ we have a natural
identification  of $V_{x,s}$ with
$V_{x,s+j \cdot \ell}$. With respect to this 
identification, we
define
$$\overline{\gg}_x := \bigoplus_{s \in \R/ {\ell \cdot \Z}} V_{x,s}.$$
Note that $\dim_\ff \overline{\gg}_x = \dim _k \gg$.
We define a product operation on $\overline{\gg}_x$ in the following
manner.
If $\overline{X}_s \in V_{x,s}$ and $\overline{X}_t
\in V_{x,t}$, then we define
$[\overline{X}_s, \overline{X}_t]$
to be the image of 
$[X_s, X_t] \in \gg_{x,(s+t)}$
 in $V_{x,(s+t)}$ where $X_s \in \gg_{x,s}$ and $X_t \in \gg_{x,t}$
are any lifts of $\overline{X}_t$ and $\overline{X}_s$, respectively. 
 Linearly extend this  operation to an operation on
$\overline{\gg}_x$.  With this product $\overline{\gg}_x$ is an
$\ff$-Lie algebra.
 For $v \in
\overline{\gg}_x$, define ${\dad}(v) \in
\End _\ff \bigl(\overline{\gg}_x \bigr)$ by ${\dad}(v)w =
[v,w]$ for all $w \in \overline{\gg}_x$. 

\vglue2pt

For more information about Hypothesis~\ref{hyp:pforallnew},
see Appendix~A.

\specialnumber{4.2.1}\proclaim{Hypothesis}  \label{hyp:pforallnew}
Suppose $x \in \BB(G)${\rm .} If $X \in \NN \cap (\gg_{x,r}
\smallsetminus 
\gg_{x,r^+})${\rm ,} then there exist $H \in \gg_{x,0}$ and $Y \in
\gg_{x,-r}$ such that  
\begin{eqnarray*}
\left[H,X\right] &=& 2X  \hbox{ mod } \gg_{x,r^+}\\
\left[H,Y\right] &=& -2Y \hbox{ mod } \gg_{x,(-r)^+}\\
\left[X,Y\right] &=& H   \hbox{ mod } \gg_{x,0^+}.
\end{eqnarray*}
If $(f,h,e)$ denotes the image of $(Y,H,X)$ in 
$V_{x,-r} \times V_{x,0} \times V_{x,r} \subset \overline{\gg}_x${\rm ,}
then $(f,h,e)$ is an ${\rmsl}_2(\ff)$\/{\rm -}\/triple{\rm ,} and 
$\overline{\gg}_x$ decomposes into a direct sum of irreducible
$(f,h,e)$\/{\rm -}\/modules of 
highest weight at most $(p-3)${\rm .}  Moreover{\rm ,}  there exists
$\bar{\lambda} \in X_*^{\ff} (\bfG_x)${\rm ,} uniquely determined up to an
element of $\X_*(\bfZ_x)$ whose differential is zero{\rm ,}  such 
that 
the following two conditions hold{\rm .}
\begin{itemize}
\ritem{1.}
The image of $d \bar{\lambda} $ in $\Lie(\bfG_x)$ coincides with the
one\/{\rm -}\/dimensional subspace spanned by $h${\rm .}
\ritem{2.} Suppose $i \in \Z${\rm .}   For   $v \in
\overline{\gg}_x$ 
$$\hbox{if $\lsup{\bar{\lambda}(t)}v = t^iv$, then $\dabs{i} \leq
(p-3)$ and $\dad(h)v = iv$.}$$  
\end{itemize}
\endproclaim

{\it Definition} 4.2.2.
In the notation of Hypothesis~\ref{hyp:pforallnew}, 
we say
that $\bar{\lambda} 
\in \X_*^{\ff}(\bfG_x)$ is {\it adapted} to the
${\rmsl}_2(\ff)$-triple obtained from the image of $(Y,H,X)$ in 
$V_{x,-r} \times V_{x,0} \times V_{x,r}$.

\specialnumber{4.2.3}\proclaim{Hypothesis} \label{hyp:decompnew}
If $X \in \NN${\rm ,} then there exists $m \in \N$ with $m \leq (p -
2)$ such that 
$\dad(X)^{m} = 0${\rm .}
\endproclaim

For more
background on the next hypothesis see, for
example,~\cite[\S5.5]{carter:finite}. \pagebreak 

\specialnumber{4.2.4}\proclaim{Hypothesis} \label{hyp:exptexist} 
Choose $m \in \N$ such that $\dad(X)^m = 0$ for all $X \in \NN${\rm .}
Suppose either that $k$ has characteristic zero or that the characteristic of
$k$ is greater than $m${\rm .}  There exists a unique
$G$\/{\rm -}\/invariant map $\exp_t \colon \NN \rightarrow \UU$ 
defined over $k$ such that for all $X \in \NN$ the adjoint action of 
$\exp_t(X)$ on $\gg$ is given by 
$$\sum_{i=0}^{m} \frac{\bigl(\dad(X)\bigr)^i}{i!}.$$
\endproclaim

 For more information about Hypothesis~\ref{hyp:morosowfork}
see~\cite[\S5.5]{carter:finite}.  

\def\scrs#1{\scriptstyle{#1}}
\specialnumber{4.2.5}\proclaim{Hypothesis} \label{hyp:morosowfork}
Suppose Hypothesis~{\rm \ref{hyp:exptexist}} is valid{\rm .}   
Suppose $X \in \NN${\rm .}
There exists an ${\rmsl}_2(k)$\/{\rm -}\/triple completing $X${\rm .}  For any
${\rmsl}_2(k)$\/{\rm -}\/triple $(Y,H,X)$ completing $X$ there is  a group
homomorphism $\varphi \colon {\bf SL}_2 \rightarrow \bG$ defined
over $k$ such that $d\varphi {\big(
\begin{array}{cc}\scrs{0}\hskip-5pt&\scrs{1}\\[-6pt]\scrs{0}\hskip-5pt&\scrs{0}\end{array} \big)} = X${\rm ,} $d\varphi
{\big( \begin{array}{cc}\scrs{0}\hskip-5pt&\scrs{0}\\[-6pt]\scrs{1}\hskip-5pt&\scrs{0}\end{array} \big)} = Y${\rm ,} 
$d\varphi {\big( \begin{array}{cc}\scrs{1}\hskip-5pt&\scrs{0}\\[-6pt]\scrs{0}\hskip-5pt&\scrs{-1}\end{array} \big)} =
H${\rm ,} and for all 
$t \in k$ 

{\rm 1.} 
$\varphi{\big( \begin{array}{cc}\scrs{1}\hskip-5pt&\scrs{t}\\[-6pt]\scrs{0}\hskip-5pt&\scrs{1}\end{array} \big)} = \exp_t( 
t X )$ and

{\rm 2.} 
$\varphi {\big( \begin{array}{cc}\scrs{1}\hskip-5pt&\scrs{0}\\[-6pt]\scrs{t}\hskip-5pt&\scrs{1}\end{array} \big)} =
\exp_t(  tY )${\rm .}
\vglue3pt\noindent
Finally{\rm ,} any two ${\rmsl}_2(k)$-triples completing $X$ are conjugate by
an element of $C_G(X)${\rm .}
\endproclaim 

{\it {R}emark} 4.2.6.
We note that the map $\varphi$ occurring in
Hypothesis~\ref{hyp:morosowfork} is uniquely determined by $d \varphi$.
\vglue8pt

For more information about Hypothesis~\ref{hyp:phimap},
see~\cite[\S1.6]{adler:thesis}. 

\specialnumber{4.2.7}\proclaim{Hypothesis} \label{hyp:phimap}
Suppose $x \in \BB(G)${\rm .}  For all  $s \in \R_{>0}$ and for all $t \in
\R$ there exists 
a map $\phi_x 
\colon \gg_{x,s} \rightarrow G_{x,s}$ such that for $V \in \gg_{x,s}$
and $W \in \gg_{x,t}$ we have
$$\lsup{\phi_x(V)}W = W + [V,W] \hbox{\, mod \,} \gg_{x,(s+t)^+}.$$ 
\endproclaim

\vglue-18pt
4.3. {\it  From Jacobson-Morosov triples over $\ff$ to
Jacobson-Morosov triples over $k$}. Fix $x \in \BB(G)$.  Suppose that $(f,h,e) \in V_{x,-r}
\times V_{x,0} \times V_{x,r} \subset \overline{\gg}_x$ is a
(nontrivial) ${\rmsl}_2(\ff)$-triple 
with adapted  $\bar{\mu} \in \X_*^{\ff}(\bfG_x)$. 
We now show that, subject to some conditions on $k$ and $\bG$, there exist $Y \in
\gg_{x,-r}$, $H \in \gg_{x,0}$ and  $X \in \gg_{x,r}$ such that
$(Y,H,X)$ is an ${\rmsl}_2(k)$-triple in $\gg$ and $(Y,H,X)$ is
a lift of $(f,h,e)$ in the obvious sense.  We 
follow~\cite[\S\S3.8--3.9]{barbasch-moy:local} where the proof is
carried out for certain $\bG$ when $r = 0$.

Let $\bS$ be a maximal $k$-split torus of $\bG$ such that $x \in
\AA(\bS,k)$.   Let $\bfS$ be the maximal 
$\ff$-split torus of $\bfG_x$ corresponding to $\bS$.  Since
maximal $\ff$-split tori are $\bfG_x(\ff)$-conjugate, there exist
$\bar{\lambda} \in \X_*(\bfS)$ and $\bar{g} \in \bfG_x(\ff)$ such that
$\bar{\lambda} = \lsup{\bar{g}}\bar{\mu}$.  Let $\lambda \in
\X_* (\bS)$ be the lift of $\bar{\lambda}$ and replace $(f,h,e)$ with
$(\lsup{\bar{g}}f, \lsup{\bar{g}}h, \lsup{\bar{g}}e)$.

For $i \in \Z$, define
$$\gg(i) := \{Z \in \gg \, | \, \lsup{\lambda (t)}Z = t^i \cdot Z 
\}  \, \hbox{ and } \,
\overline{\gg}_x(i) := \{v \in \overline{\gg}_x \, | \,
\lsup{\bar{\lambda} (t)}v = t^i \cdot v 
\}.$$
For $i \in \Z$ and  $s \in \R$ define
$$\gg_{x,s}(i) := \{Z \in \gg_{x,s} \, | \, \lsup{\lambda (t)}Z = t^i
\cdot Z 
\}
\, \hbox{ and } \,
V_{x,s}(i) := \{v \in V_{x,s} \, | \, \lsup{\bar{\lambda} (t)}v = t^i
\cdot v 
\}.$$
Because $x \in \AA(\bS,k)$, $\lambda \in \X_*(\bS)$, and 
$\bar{\lambda} \in \X_*(\bfS)$ we have
$$\gg_{x,s} = \bigoplus_i \gg_{x,s}(i)  \hbox{\, and \,} V_{x,s} =
\bigoplus_i V_{x,s}(i)$$
for $s \in \R$ and $i \in \Z$.

\specialnumber{4.3.1}\proclaim{Lemma} \label{lem:Xsquareiso}
Suppose that Hypothesis~{\rm \ref{hyp:pforallnew}} holds{\rm .}
If $X \in \gg_{x,r}(2)$ is any lift of $e${\rm ,} then  for all $s \in \R${\rm ,} the map
$$\dad(X)^2 \colon \gg_{x,s-r}(-2) \rightarrow \gg_{x,s+r}(2)$$
is an isomorphism{\rm .} 
\endproclaim

{\it Proof}.
Fix $X \in \gg_{x,r}(2)$ which is a lift of $e$. Note that $X$ is nilpotent. 
Since $\dad(X)^2$ takes $\gg(-2)$ to $\gg(2)$ and $k$ is complete, it
will be sufficient to show that for all $t \in \R$, the map
$$\dad(e)^2 \colon V_{x,t-r}(-2) \rightarrow V_{x,t+r}(2)$$
is an isomorphism.

From Hypothesis~\ref{hyp:pforallnew} we have
that the space $\overline{\gg}_x$  is a direct sum of irreducible
$(f,h,e)$-modules. 
Consequently, it follows from ${\rmsl}_2(\ff)$-representation
theory that the map
$$\dad(e)^2 \colon \overline{\gg}_{x}(-2) \rightarrow \overline{\gg}_{x}(2)$$
is an isomorphism.  The result follows.
\hfill\qed

\specialnumber{4.3.2}\proclaim{{C}orollary} \label{cor:liftsexist}
Suppose that Hypotheses~{\rm \ref{hyp:pforallnew}} and~{\rm \ref{hyp:decompnew}}
hold{\rm .}  If $X \in 
\gg_{x,r}(2)$ is a lift of $e${\rm ,} then there exist  lifts
$Y \in \gg_{x,-r}$ of $f$ and  $H \in \gg_{x,0}$ of $h$ such that
$(Y,H,X)$ is an 
${\rmsl}_2(k)$\/{\rm -}\/triple in $\gg${\rm .}  
\endproclaim

{\it Proof}.
Fix $X \in \gg_{x,r}(2)$ which is a lift of $e$.
Since $\dad(X)^2 \colon \gg_{x,-r}(-2) \rightarrow \gg_{x,r}(2)$ is a
surjection
there exists $Y \in \gg_{x,-r}(-2)$
such that
$\dad(X)^2 Y = -2X$.
Since  $\dad(e)^2 \colon \overline{\gg}_{x}(-2) \rightarrow
\overline{\gg}_{x}(2)$  is
injective and $\dad(e)^2 f = -2e$, 
$Y$ is necessarily a lift of $f$.
Let $H = [X,Y] \in \gg_{x,0}(0)$. $H$ is a lift of $h$.  We
also have $[H,X] = 2X$.  We need to check that $[H,Y] = -2Y$.  

It follows from Hypothesis~\ref{hyp:decompnew} and
Morosov's 
theorem (see~\cite[Lemma~7, p.~98]{jacobson:lie}
or~\cite[the proof of Proposition~5.3.1; in particular, 
pp.~140--141]{carter:finite}) that there exists $Y' \in \gg$ such that 
$(Y',H,X)$ is an ${\rmsl}_2(k)$-triple.  Since $H \in \gg(0)$,
$X \in \gg(2)$, and $[\gg(i), \gg(j)] \subset \gg(i+j)$, we can assume
that $Y' \in \gg(-2)$.  However, since $\dad(X)^2 \colon
\gg_{x,s-r}(-2) \rightarrow \gg_{x,s+r}(2)$ is injective for all $s \in
\R$, we must have $Y' = Y$.
\hfill\qed\pagebreak

4.4. {\it Some fixed-point results for one-parameter subgroups}.
Fix $\lambda \in \X_*^k(\bG)$.  Let $\bM$ denote the $k$-{L}evi
subgroup of $\bG$ whose group of 
$k$-rational points is the Levi subgroup $M = C_\Gn(\lambda)$.   
Note that for all $z \in \BB(M)$ we have $d \lambda (R^\times) \subset 
\gg_{z,0} \smallsetminus \gg_{z,0^+}$. 

The following lemma and its subsequent applications in
Corollaries~\ref{cor:centisfixn} and~\ref{cor:liecentisfixn} 
arose from discussions with Gopal Prasad.  The results of this subsection
may be thought of as natural generalizations of the material
in\break \cite[\S3.6]{corvallis}. 

\specialnumber{4.4.1}\proclaim{Lemma} \label{lem:twoorall}
Suppose $F$ is a
codimension one $0$\/{\rm -}\/facet in $\BB(\bM,K)${\rm .}  Either every
$0$\/{\rm -}\/alcove of $\BB(\bG,K)$ which contains $F$ in its closure belongs to
$\BB(\bM,K)$ or exactly two of the $0$\/{\rm -}\/alcoves which contain $F$ in
their closure lie in $\BB(\bM,K)${\rm .}
\endproclaim

{\it Proof}.
Since $F 
\subset \BB(\bM,K)$, there exists a maximal $K$-split torus $\bT$ such
that $\bT \subset \bM$ and $F \subset \AA(\bT,K)$.  Since $\lambda \in
\X_*(\bT)$, we can consider the image $\bar{\lambda}$ of $\lambda$ in
$\X_*(\bfG_F)$.   Let $C_1$ and $C_2$ denote the $0$-alcoves in
$\AA(\bT,K)$ which contain $F$ in their closure.  Let $\bfB$ denote the
Borel subgroup of $\bfG_F$ corresponding to~$C_1$.

First suppose that $\bar{\lambda}$ lies in the center of
$\bfG_F$.  Since every Borel subgroup of $\bfG_F$ is conjugate to
$\bfB$ by an element of $\bfG_F$, we conclude that every $0$-alcove in
$\BB(\bG,K)$ which contains $F$ in its closure is conjugate to $C_1$
by an element of $\bG(K)_F \cap \bM(K)$.

Now suppose that $\bar{\lambda}$ does not lie in the center of $\bfG_F$.
Since the derived group
of $\bfG_F$ is either ${\bf SL}_2$ or ${\bf PGL}_2$, 
it follows
that there are 
exactly two Borel subgroups in $\bfG_F$ containing $\bar{\lambda}$.
These Borel subgroups correspond to $C_1$ and~$C_2$.
\hfill\qed

\vglue9pt
\specialnumber{4.4.2}\proclaim{{C}orollary} \label{cor:centisfixn}
Suppose that $\ff$ has more than three elements{\rm .}  We have  
$$\BB(M) = \BB(G)^{\lambda(R^\times)}.$$ 
\endproclaim

{\it Proof}.
``$\subset$'':  
Since any natural embedding of $\BB(M)$ into $\BB(G)$ is
$M$-equivariant, this follows from the fact that $\lambda(R^\times)$
fixes $\BB(M)$ point-wise. 

``$\supset$'': Suppose $\BB(M) \subsetneq
\BB(G)^{\lambda(R^\times)}$.  We will obtain a contradiction.

Since  
the group $\lambda(R^\times)$
fixes $\BB(\bM,K)$, there exists a $0$-alcove
$C$ in $\BB(\bG,K)$ such that $\bar{C} \cap \BB(\bM,K)$ has codimension one and
$\lambda(R^\times)$ fixes $\bar{C}$.

Choose an apartment $\AA$ in  $\BB(\bM,K)$ such that
$F = \bar{C} \cap \BB(\bM,K) \subset \AA$.  Let $\bT$ be the maximal
$K$-split torus in $\bG$ corresponding to $\AA$.  Note that $\lambda
\in \X_*(\bT)$.  Since $C$ does not lie in
$\BB(\bM,K)$, from the proof of Lemma~\ref{lem:twoorall} we conclude
that the image of $\lambda$ in $\bfG_F$ does not lie in the center of $\bfG_F$.
Thus, since the derived group \pagebreak of $\bfG_F$ is
either ${\bf SL}_2$ or ${\bf PGL}_2$, we conclude that if the
cardinality of $\ff$ is greater than three, then the image of
$\lambda(R^\times)$ in 
$\bfG_F$ lies only in those Borel subgroups of $\bfG_F$ corresponding to
$0$-alcoves in $\BB(\bM,K)$.  That is, with our restrictions on $\ff$,
$C$ cannot be fixed by $\lambda(R^\times)$.
\hfill\qed
\vglue12pt
\specialnumber{4.4.3}\proclaim{{C}orollary} \label{cor:liecentisfixn}
Suppose the characteristic of $\ff$ is not two{\rm .}  Let $H = d
\lambda (1)${\rm .} 
If $y \in \BB(G)${\rm ,} then
$H \in
\gg_{y,0}$ if and only if $y \in \BB(M)${\rm .}
\endproclaim
\vglue12pt
{\it Proof}.
``$\Leftarrow$'': This is immediate.

``$\Rightarrow$'': Let ${\cal  C}$ denote the set of $z \in
\BB(\bG,K)$ for which 
$H \in {\rmg }(K)_{z,0}$.
The set ${\cal  C}$ is convex and contains $\BB(\bM,K)$.

Suppose that ${\cal  C} \neq \BB(\bM,K)$.  We will derive a contradiction.
Since ${\cal  C}$ is convex and contains $\BB(\bM,K)$, there exists
a $0$-alcove 
$C$ in $\BB(\bG,K)$ such that $\bar{C} \cap \BB(\bM,K)$ has codimension one and
$H \in \rmg(K)_{c,0}$ for all $c \in C$.

Choose an apartment $\AA$ in $\BB(\bM,K)$ such that
$\bar{C} \cap \BB(\bM,K) \subset \AA$.  Let $\bT$ be the maximal $K$-split
torus in $\bG$ corresponding to $\AA$.  Note that $\lambda \in
\X_*(\bT)$.  
Since $C$ does not lie in
$\BB(\bM,K)$, from the proof of Lemma~\ref{lem:twoorall} we conclude
that the image of $\lambda$ in $\bfG_F$ does not lie in the center of $\bfG_F$.
Thus, since the derived group of $\bfG_F$ is
either ${\bf SL}_2$ or ${\bf PGL}_2$ and since the
characteristic of $\ff$ is not two, the image of
$H$ in 
$\Lie(\bfG_F)$ lies only in those Borel subalgebras of $\Lie(\bfG_F)$
corresponding to 
$0$-alcoves in $\BB(\bM,K)$.  That is, with our restrictions on $\ff$,
we cannot have $H \in {\rmg }(K)_{c,o}$ for any $c \in C$.
\hfill\qed
\vglue12pt

4.5. {\it From Jacobson-Morosov triples over $k$ to Jacobson-Morosov
triples over $\ff$}.  Different versions of Lemma~\ref{lem:kushnirsky} were proved
independently (and nearly  
simultaneously) by Eugene Kushnirsky and myself. 
The proof here is
due to Eugene Kushnirsky;  I thank him for allowing me to publish it
here.  My proof will appear elsewhere.

\specialnumber{4.5.1}\proclaim{Lemma} \label{lem:kushnirsky}
For any $x,y \in \BB(G)${\rm ,} we have $\stab_G(x) \cap G_y = G_x \cap
\stab_G (y) = G_{\{x,y\}}${\rm .}
\endproclaim

Here $G_{\{x,y\}}$ is the $\Gal(K/k)$-fixed points of the group of
$R_K$-rational points of the identity component of the group scheme
associated to the set 
$\{x,y\}$ (see~\cite[\S3.4]{corvallis} and~\cite[\S1.2.12]{bruhat-tits:two}).
\vglue12pt

{\it Proof {\rm (}\/Eugene Kushnirsky\/{\rm )}}.
Without loss of generality, we work over~$K$.  Let $\AA$ be
an apartment in $\BB(\bG,K)$ containing $x$ and 
$y$.  Let $\bT$ denote the maximal $K$-split torus of $\bG$
corresponding to $\AA$ and let $\bZ$ denote the centralizer in $\bGn$ of $\bT$.

For $\alpha \in \Phi(\bT,K)$, let $U_\alpha \subset \bGn(K)$ denote the
corresponding root subgroup.
For a fixed ordering on $\Phi(\bT,K)$ we define $U^+$ (resp., $U^-$)
to be the group generated by $\{U_\alpha \}_{\alpha > 0}$ (resp.,
$\{U_\alpha \}_{\alpha < 0}$).  For $z \in \AA$ we define $U_z^{\pm} =
\stab_G(z) \cap U^\pm$.   
Choose an ordering on $\Phi(\bT, K)$ so that $U_y^+ \subset U_x^+$;
this implies that $U^+_{\{x,y\}} = U_y^+$.  According
to~\cite[{\it Corollaire}~4.6.7]{bruhat-tits:two}, $G_y = U_y^+U_y^-U_y^+
{\cal  I}^0(R_K)$.  Here ${\cal  I}$ is the smooth $R_K$-model for
$\bZ$ constructed in~\cite[\S4.4]{bruhat-tits:two}.  Another application
of~\cite[{\it Corollaire}~4.6.7]{bruhat-tits:two} produces 
\begin{eqnarray*} \noalign{\vskip-4pt}
\stab_G(x) \cap G_y &=& (\stab_G(x) \cap U_y^+U_y^-U_y^+){\cal  I}^0(R_K)\\[-2pt]
&=& U_{\{x,y\}}^+(\stab_G(x) \cap U_y^-) U_{\{x,y\}}^+
{\cal  I}^0(R_K)\\[-2pt]
&=& U^+_{\{x,y\}}U^-_{\{x,y\}}U_{\{x,y\}}^+ {\cal  I}^0(R_K)\\[-2pt]
&=& G_{\{x,y\}}. \\ \noalign{\vskip-36pt}
\end{eqnarray*}
\hfill\qed 

\specialnumber{4.5.2}\numbereddemo{{R}emark}  \label{rem:kushnirsky}
Recall the definition of unipotent in Section~2.1.
Since every unipotent element belongs to some parahoric subgroup, Lemma~\ref{lem:kushnirsky} implies that if $u \in G$ is unipotent and
$x \in \BB(G)$ such that $ux = x$, then $u \in G_x$.
\enddemo

For the remainder of this subsection we fix a nontrivial $X \in
\NN$, and  we suppose that
Hypothesis~\ref{hyp:morosowfork} holds. Let $(Y,H,X)$ be an
${\rmsl}_2(k)$-triple completing $X$.  Suppose that $\varphi$ is 
a homomorphism for $(Y,H,X)$
as described in Hypothesis~\ref{hyp:morosowfork}.  We have $H = d \varphi
({\big(
\begin{array}{cc}\scrs{1}\hskip-5pt&\scrs{0}\\[-6pt]\scrs{0}\hskip-5pt&\scrs{-1}\end{array} \big)})$ and $Y = d \varphi
({\big(
\begin{array}{cc}\scrs{0}\hskip-5pt&\scrs{0}\\[-6pt]\scrs{1}\hskip-5pt&\scrs{0}\end{array} \big)})$.   We wish to find
a point $y \in \BB(G)$ such that $Y \in \gg_{y,-r}$, $H \in \gg_{y,0}$,
and $X \in \gg_{y,r}$.  

I thank Gopal Prasad for explaining to me the proof of the following
lemma; this lemma occurs without proof
in~\cite[Corollary~3.7~(1)]{barbasch-moy:local}.

\proclaimtitle{Dan Barbasch and Allen Moy}
\specialnumber{4.5.3}\proclaim{Lemma}\label{lem:mpalmost}
Suppose Hypothesis~{\rm \ref{hyp:morosowfork}} holds{\rm .}  There exists $x \in
\BB(G)$ such that $Y,H,X \in \gg_{x,0}${\rm .}
\endproclaim

{\it {P}roof {\rm (}\/Gopal Prasad\/{\rm )}}.
Let $J = \varphi \bigl({\bf SL}_2(R_K)\bigr) \subset \bGn(K)$.  The 
group $J \rtimes \Gal(K/k)$ acts on $\BB(\bG,K)$.  Moreover, since 
 $J \rtimes \Gal(K/k)$ is bounded, its action has a fixed
point~\cite[\S2.3.1]{corvallis}.  Let $x \in \BB(\bG,K)$ be such a fixed-point.

Let ${\cal  G}$ denote the $R$-group scheme associated to
$\stab_{\bGn(K)}(x)$ 
(see~\cite{bruhat-tits:two}).  The generic fiber ${\cal  G}
\otimes_R k$ is $\bGn$ and the group of $R_K$-rational points of
${\cal  G}$ is $\stab_{\bGn(K)}(x)$.  Let $L({\cal  G})$ denote the
Lie algebra of ${\cal  G}$.  $L({\cal  G})$ is a lattice in $\gg$
and ${\rmg }(K)_{x,0} = L({\cal  G}) \otimes_R R_K$.  Let ${\cal  J}$
denote the $R$-group scheme associated to ${\bf SL_2}(R_K)$.  
From~\cite[{Proposition}~1.7.6]{bruhat-tits:two} the map
$\varphi$ induces an $R_K$-scheme homomorphism of ${\cal  J}$ into
${\cal  G}$.  
 Consequently, $d \varphi \bigl({\rmsl}_2(R_K)\bigr) \subset
{\rmg }(K)_{x,0}$.   

Since $x$ is fixed by $\Gal(K/k)$, we have $x \in \BB(G)$
and $Y,H,X \in \gg_{x,0}$.\hfill\qed

{\it {R}emark} 4.5.4.
The images of $Y$, $H$, and $X$ in $V_{x,0}$ form an
${\rmsl}_2(\ff)$-triple.
 
\specialnumber{4.5.5}\proclaim{{C}orollary}  \label{cor:mpalmostbeefnew}
Suppose Hypotheses~{\rm \ref{hyp:decompnew}} and~{\rm \ref{hyp:morosowfork}} hold{\rm .}  
If $x \in \BB(G) = \BB(\bG,K)^{\Gal(K/k)}${\rm ,} then
$$x \in \BB(\bG,K)^{\varphi({\bf SL}_2(R_K))} 
\hbox{\, if and
only if \,} d\varphi\bigr({\rmsl}_2(R)\bigl) \subset \gg_{x,0}.$$
\endproclaim

{\it Proof}.
``$\Rightarrow$'':   This follows from the proof of
Lemma~\ref{lem:mpalmost}.

``$\Leftarrow$'': Suppose $F$ is a $\Gal(K/k)$-invariant $0$-facet in $\BB(\bG,K)$ such that $x
\in F^{\Gal(K/k)}$.  
Since $d \varphi \bigl( {\rmsl}_2(R) \bigr) \subset \gg_{x,0}$,
it follows from Hypotheses~\ref{hyp:decompnew} and~\ref{hyp:morosowfork}
that for all $t \in 
R_K$ both $\varphi ({\big(
\begin{array}{cc}\scrs{1}\hskip-5pt&\scrs{t}\\[-6pt]\scrs{0}\hskip-5pt&\scrs{1}\end{array} \big)})$ and 
$\varphi ( {\big(
\begin{array}{cc}\scrs{1}\hskip-5pt&\scrs{0}\\[-6pt]\scrs{t}\hskip-5pt&\scrs{1}\end{array} \big)})$ lie in
$N_{\bG(K)}({\rmg }(K)_{x,0}) \cap N_{\bG(K)}({\rmg }(K)_{x,0^+})$.  Since these elements generate
$\varphi\bigl({\bf SL}_2(R_K)\bigr)$, 
it follows from  Lemma~\ref{lem:normstab} that
$\varphi\bigl({\bf SL}_2(R_K)\bigr) 
\subset \stab_{\bGn(K)}(F)$.  From Remark~\ref{rem:kushnirsky} it follows that
$\varphi\bigl({\bf SL}_2(R_K)\bigr) \subset \stab_{\bGn(K)}(F)$
if and only if
$\varphi\bigl({\bf SL}_2(R_K)\bigr) \subset \bG(K)_F$.  Thus,
$\varphi\bigl({\bf SL}_2(R_K)\bigr) \subset \bG(K)_F =
\bG(K)_x \subset \stab_{\bG(K)}(x).$
\hfill\qed\vglue8pt

Let $\lambda
\in \X_*^k(\bG)$ be the one-parameter subgroup derived from $\varphi$.
That is, $\lambda(t) = \varphi ({\big(
\begin{array}{cc}\scrs{t}\hskip-5pt&\scrs{0}\\[-6pt]\scrs{0}\hskip-5pt&\scrs{t^{-1}}\end{array} \big)})$ for all $t \in
k^\times$.

\specialnumber{4.5.6}\numbereddemo{Definition}\label{defn:fadapt}
The one-parameter subgroup $\lambda$ constructed in the preceding
paragraph is said to be {\it adapted} to the
${\rmsl}_2(k)$-triple $(Y,H,X)$.
\enddemo

Define the Levi subgroup $M = C_\Gn(\lambda)$.    We now present two
corollaries of the results in subsection~4.4.

\specialnumber{4.5.7}\proclaim{{C}orollary} \label{cor:ftoffn}
Suppose Hypotheses~{\rm \ref{hyp:decompnew}} and~{\rm \ref{hyp:morosowfork}} hold{\rm .} 
There exists $y \in \BB(G)$ such that $Y \in \gg_{y,-r}${\rm ,} $H \in
\gg_{y,0}${\rm ,} and $X \in \gg_{y,r}${\rm .}
\endproclaim

\demo{Proof}
The hypotheses imply that $\ff$ has more than three elements.

Choose $x \in \BB(G)$ as in Lemma~\ref{lem:mpalmost}.  Since $x$ is
fixed by the group $\lambda(R^\times)$, Corollary~\ref{cor:centisfixn}
tells us $x \in \BB(M)$.  Fix an apartment $\AA$ in
$\BB(M)$ which contains $x$.  Since the group $\lambda(k^\times)$ lies in the center of
$M$, $\lambda$ acts on every apartment in $\BB(M)$ by translation.
It therefore makes sense to define $y = x + (r/2) \cdot \lambda \in
\AA$.  This $y$ 
satisfies the requirements of the lemma.
\enddemo

\specialnumber{4.5.8}\numbereddemo{{R}emark}
The image of $(Y,H,X)$ in $V_{y,-r} \times V_{y,0} \times V_{y,r}$ forms an
${\rmsl}_2(\ff)$-triple under the inherited Lie algebra operation.
\enddemo

The following result
may be interpreted as a sharpening of Corollary~\ref{cor:ftoffn}.

\specialnumber{4.5.9}\proclaim{{C}orollary}  \label{cor:liecentisfixm}
Suppose Hypotheses~{\rm \ref{hyp:decompnew}} and~{\rm \ref{hyp:morosowfork}} hold{\rm .} 
If $y \in \BB(G)$ such that 
$Y \in \gg_{y,-r}${\rm ,} $H \in
\gg_{y,0}${\rm ,} and $X \in \gg_{y,r}${\rm ,} then $y$ must lie in $\BB(M)${\rm .}
\endproclaim

\demo{Proof}
The hypotheses imply that the characteristic of $\ff$ is not two.
Therefore, the result follows from Corollary~\ref{cor:liecentisfixn}.
\enddemo\pagebreak

\section{The parametrization}  \label{sec:parametrize}

Fix $r \in \R$.
In this section we combine the material of the previous two sections
and produce a parametrization of the nilpotent orbits in $\gg$.

\demo{{\rm 5.1.} A {\rm ``}\/building set\/{\rm ''} related to an
${\rmsl}_2(k)$-triple}
Suppose Hypotheses~\ref{hyp:decompnew} and ~\ref{hyp:morosowfork} hold.
Given an ${\rmsl}_2(k)$-triple in $\gg$, we want to produce a
nice subset of $\BB(G)$. The idea for the definitions in this subsection
originated in~\cite{moy-prasad:refined}.  

Fix $Z \in \NN$ and $s \in \R$.

\specialnumber{5.1.1}\numbereddemo{Definition}
$$ \BB(Z,s) := \{z \in \BB(G) \, | \, Z \in \gg_{z,s} \}.$$
\enddemo

The set $\BB(Z,s)$ is a nonempty and convex subset of $\BB(G)$.
Moreover, it is the union of generalized $s$-facets. 

\specialnumber{5.1.2}\proclaim{Lemma}  \label{lem:buildclose}
$\BB(Z,s)$ is closed{\rm .}
\endproclaim

\demo{Proof}
Suppose that $y \in \overline{\BB(Z,s)}$.
Let $F^*$ be the generalized $s$-facet containing $y$.  From
Lemma~\ref{lem:opengen},  the
union of all generalized $s$-facets that contain $y$ in their
closure is an open neighborhood of $y$.
 Consequently, there
exists a generalized $s$-facet $H^*$ such that $H^* \subset \BB(Z,s)$
and $F^* \subset \overline{H^*}$.  From Corollary~\ref{cor:closure},
we have $Z \in \gg_{H^*,s} \subset \gg_{F^*,s} = \gg_{y,s}$.  Thus $y
\in \BB(Z,s)$.
\enddemo

Suppose $(Y,H,X)$ is an (possibly trivial) ${\rmsl}_2(k)$-triple in
$\gg$. 

\specialnumber{5.1.3}\numbereddemo{Definition}
$$\BB(Y,H,X) := \BB(X,r) \cap \BB(Y,-r).$$
\enddemo

The set $\BB(Y,H,X)$ is a nonempty (Corollary~\ref{cor:ftoffn}), closed
(Lemma~\ref{lem:buildclose}), and convex subset of
$\BB(G)$. Moreover, it is also the 
union of generalized\break $r$-facets.  

These properties imply the following
result (see also~\cite[\protect{Lemma~3.6}]{barbasch-moy:local}).

\specialnumber{5.1.4}\proclaim{Lemma} \label{lem:maximalgenr}
Suppose $F_1^*, F_2^* \in {\cal  F}(r)$ and $F_i^* \subset
\BB(Y,H,X)${\rm .}  If $F_1^*$ and $F_2^*$ are maximal
generalized $r$\/{\rm -}\/facets in $\BB(Y,H,X)${\rm ,} then $F_1^*$ and $F_2^*$ are
strongly $r$\/{\rm -}\/associated{\rm .}
\endproclaim

\demo{Proof}
Choose $x_i \in F_i^*$.  Let $\AA$ be an apartment in $\BB(G)$
containing $x_1$ and $x_2$.  Since $F^*_i$ is maximal in $\BB(Y,H,X)$
and since $\BB(Y,H,X)$ is convex, we have $\emptyset \neq F^*_j \cap \AA
\subset A(\AA, F^*_i)$ for $i,j \in \{1,2\}$.
\enddemo

\specialnumber{5.1.5}\numbereddemo{{R}emark}
In the language of Section~4.5, if $X$ is not trivial, then from
Corollaries~\ref{cor:centisfixn},~\ref{cor:mpalmostbeefnew}
and~\ref{cor:liecentisfixm}  we have 
$$\BB(Y,H,X) = \BB(\bG,K)^{\varphi({\bf SL}_2(R_K)) \rtimes
\Gal(K/k)} + (r/2)\cdot \lambda.$$ 
The sum on 
the right-hand side occurs in $\BB(G)^{\lambda(R^\times)} = \BB\bigl(
C_{\Gn}(\lambda)\bigr)$. 
\enddemo

5.2. {\it An extension of some work of J.\/{\rm -}\/L.~Waldspurger}.
We assume that Hypotheses~\ref{hyp:decompnew},~\ref{hyp:morosowfork},
and~\ref{hyp:phimap}  are in effect. 
\advance\eqcount by 1

Fix $X \in \NN \smallsetminus \{0\}$.  
Suppose $(Y,H,X)$ is an ${\rmsl}_2(k)$-triple in $\gg$. Suppose
$\lambda \in \X_*^k(\bG)$ is adapted to $(Y,H,X)$. 
Fix $x \in \BB(Y,H,X)$. 

We will explore the relationship between the coset $X +
\gg_{x,r^+}$ and the nilpotent orbit $\lsup{G}X$.
For example, from~\cite{barbasch-moy:local} we  
expect that $\lsup{G}X$ is the unique nilpotent orbit of minimal
dimension which intersects $X+ \gg_{x,r^+}$ nontrivially.  This result
requires some work;  
we follow 
J.-L.~Waldspurger's
presentation\break \cite[\protect{\S{IX.4}}]{waldspurger:integrales}.   

Since $H \in \gg_{x,0}$,
it follows from Corollary~\ref{cor:liecentisfixm} that there exists a 
maximal $k$-split torus $\bS$ in $\bG$ 
such that $x \in \AA(\bS,k)$ and $\lambda \in \X_*(\bS)$.
For $i \in \Z$ and $s \in \R$, define (for $\lambda$) the objects
$\gg(i)$ and  
$\gg_{x,s}(i)$ as in Section~4.3.
As before we have 
\begin{equation} \label{equ:idecomp}
\gg_{x,s} = \bigoplus_i \gg_{x,s}(i).
\end{equation}

\specialnumber{5.2.1}\proclaim{Lemma}
Assume that Hypotheses~{\rm \ref{hyp:decompnew},}
{\rm \ref{hyp:morosowfork},} and~{\rm \ref{hyp:phimap}}  hold{\rm .}
$$\lsup{G_x^+}\bigl(X + C_{\gg_{x,r^+}}(Y)\bigr) = X + \gg_{x,r^+}.$$
\endproclaim

\phantom{odd}
\vglue-36pt
{\it {P}roof {\rm (}\/A generalization of an argument of J.\/{\rm -}\/L.~Waldspurger\/{\rm )}}.

``$\subset$'':  There is nothing to prove here.

``$\supset$'':
From Hypothesis~\ref{hyp:decompnew}
and~\cite[Proposition~5.4.8]{carter:finite} we can write $\gg$ as a
direct sum of 
irreducible $(Y,H,X)$-modules of highest weight at most $(p - 3)$.  Write
$$\gg = \bigoplus_{\rho \in \Z} \gg_\rho$$
where $\gg_\rho$ denotes the isotypic component consisting
of irreducible $(Y,H,X)$-modules of
highest weight $\rho$.
For $i, \rho \in \Z$ we define $\gg(\rho,i) := \gg_\rho \cap
\gg(i)$.  We have $\gg(i) = \oplus_\rho \gg(\rho,i)$ and so $\gg =
\oplus_{\rho,i}\gg(\rho,i)$. 
For $i,\rho \in \Z$
and $s \in \R$ we define 
$$\gg_{x,s}(\rho, i) := \gg_{x,s} \cap \gg(\rho,i).$$
We first want to show
\begin{equation} \label{equ:ijdecomp}
\gg_{x,s} = \bigoplus_{\rho,i} \gg_{x,s}(\rho,i).
\end{equation}
\pagebreak

A calculation shows that if $\gg(\rho, i)$ is nontrivial, then
$$\gg(\rho,i) = \{ Z \in \gg(i) \, | \, \bigl(\dad(X) \circ
\dad(Y)\bigr) (Z) =   j(\rho,i) \cdot Z \},$$ 
where $j(\rho,i) := \bigl(
(\rho + 1)^2 - (i-1)^2\bigr)/4$.
Note that if 
$\gg(\rho,i)$  
and $\gg(\rho',i)$ are
nontrivial and $\rho \neq \rho'$, then $(j(\rho,i) - j(\rho',i),p)=1$. 

Fix $\rho, i \in \Z$ such that $\gg(\rho, i)$ is nontrivial.  Define the
nonzero integer 
$$C(\rho) := \prod_{\rho' \neq \rho; \, \, \, \gg(\rho',i) \neq \{0\}}
\bigl( j(\rho,i) - j(\rho',i)\bigr).$$ 
From the previous paragraph $(C(\rho),p) = 1$ and so $C(\rho) \in R^\times$.  
Write $Z \in \gg(i)$ as $Z = \sum Z_{\rho'}$ where $Z_{\rho'} \in
\gg(\rho',i)$.  Since the operator
$$C(\rho)\inv \cdot \prod_{\rho' \neq \rho; \, \, \, \gg(\rho',i) \neq
\{0\}} \bigl(\dad(X) \circ 
\dad(Y) -  j(\rho',i)\bigr)$$
maps $Z$ to $Z_\rho$ and
preserves depth, we have 
\begin{equation} \label{equ:jdecomp}
\gg_{x,s}(i) = \bigoplus_\rho \gg_{x,s}(\rho,i).
\end{equation}
Equations~(\ref{equ:jdecomp}) and~(\ref{equ:idecomp}) imply that
equation~(\ref{equ:ijdecomp}) is valid. 
 
Note that $C_\gg(X) = \oplus_i \gg(i,i)$ and $C_\gg(Y) = \oplus_i
\gg(-i,i)$.  Equation~(\ref{equ:ijdecomp}) 
tells us that 
\begin{equation} \label{equ:centats}
C_{\gg_{x,s}}(Y) = \bigoplus_i \gg_{x,s}(-i,i). 
\end{equation}
From equation~(\ref{equ:jdecomp}) and its proof we have 
$$\gg_{x,s}(i) = \gg_{x,s}(-i,i) + \dad(X)
\bigl(\gg_{x,(s-r)}(i-2)\bigr).
$$
Combining 
this, equation~(\ref{equ:idecomp}), and equation~(\ref{equ:centats})
yields
\begin{equation} \label{equ:final}
\gg_{x,s} = C_{\gg_{x,s}}(Y) + \dad(X) (\gg_{x,(s-r)}).
\end{equation}

Suppose $Z \in \gg_{x,r^+}$.  We wish to produce an $h \in G_x^+$ and
$C \in C_{\gg_{x,r^+}}(Y)$ such that $\lsup{h}(X+C) = X + Z$. Let $h_0 =
1$ and $C_0 = 0$. 

Fix $s_1 > r$ such that $\gg_{x,r^+} = \gg_{x,s_1} \neq
\gg_{x,s_1^+}$.  From equation~(\ref{equ:final}), we can write $Z =
C'_1 + \dad(X) P_1$ with $C'_1 \in C_{\gg_{x,s_1}}(Y)$ and $P_1 \in
\gg_{x,(s_1 - r)}$.  From Hypothesis~\ref{hyp:phimap}, there exists
$h'_1 = \phi_x(-P_1) \in G_{x,(s_1 -r)} \subset G_x^+$ such that 
\begin{eqnarray*}
\lsup{h'_1 h_0}(X + C_0 + C'_1) &=& X + C'_1 + \dad(X)P_1 \hbox{ \, mod \, }
\gg_{x,s_1^+}\\ 
&=& X + Z   - Z_1
\end{eqnarray*}
with $Z_1 \in \gg_{x,s_1^+}$. Let $h_1 = h'_1 \cdot h_0$ and $C_1 =
C_0 + C'_1$

Fix $s_2 > s_1$ such that  $\gg_{x,s_1^+} = \gg_{x,s_2} \neq
\gg_{x,s_2^+}$.  From equation~(\ref{equ:final}), we can write $Z_1 =
C'_2 + \dad(X) P_2$ with $C'_2 \in C_{\gg_{x,s_2}}(Y)$ and $P_2 \in
\gg_{x,(s_2 - r)}$.  From Hypothesis~\ref{hyp:phimap}, there exists
$h'_2 = \phi_x(-P_2) \in G_{x,(s_2-r)} \subset G_{x,(s_1 -r)}$ such that
\begin{eqnarray*}
\lsup{h'_2 h_1}(X + C_1 + C'_2) &=& \lsup{h'_2}(X + Z - Z_1  + C'_2)
\hbox{ \, mod \, } 
\gg_{x,s_2^+}\\[6pt]
&=& X + Z  - Z_1 +C'_2 + \dad(X)P_2 \hbox{\, mod \,} \gg_{x,s_2^+}\\[6pt]
&= &X + Z -Z_2
\end{eqnarray*}
with $Z_2 \in \gg_{x,s_2^+}$.  Let $h_2 = h'_2 \cdot h_1$ and $C_2 =
C_1 + C'_2$.

Continuing in this way we produce a sequence $r < s_1 < s_2
< \cdots < s_n \cdots$ with $s_n \rightarrow \infty$, elements $h_n
= h'_n h_{(n-1)} \in G^+_x$  with $h'_n \in G_{x,(s_n -r)}$
, and elements $C_n = C_{(n-1)} + C_n'  \in C_{\gg_{x,r^+}}(Y)$
with $C_n' \in C_{\gg_{x,s_n}}(Y)$ such that 
$$\lsup{h_n}(X + C_n) = X + Z
\hbox{\, mod \,} \gg_{x,s_n^+}.$$
  Let $h = \lim_{n \rightarrow \infty}h_n$ and $C
= \lim_{n \rightarrow \infty}C_n$.  Then $h \in G_x^+$, $C \in
C_{\gg_{x,r^+}}(Y)$ and $\lsup{h}(X+C) = X+ Z$.
\hfill\qed

\specialnumber{5.2.2}\proclaim{Lemma} 
Suppose that Hypothesis~{\rm \ref{hyp:morosowfork}} is valid{\rm .}
$$\bigl(X + C_{\gg}(Y)\bigr) \cap \lsup{G}X = \{X\}.$$
\endproclaim

\vglue-8pt
{\it Proof}.
See, for example,~\cite[\protect{V.7~(9)}]{waldspurger:integrales}.
\hfill\qed
 \vglue12pt

{\elevensc Corollary 5.2.3.}  
{\it Assume that Hypotheses~{\rm \ref{hyp:decompnew}, \ref{hyp:morosowfork},}
and~{\rm \ref{hyp:phimap}} hold}. 
\vglue8pt 
\hfill $(X + \gg_{x,r^+}) \cap \lsup{G}X = \lsup{G_x^+}X. $\hfill
  
\specialnumber{5.2.4}\proclaim{{C}orollary} \label{cor:uniqueorbit}
Assume that Hypotheses~{\rm \ref{hyp:decompnew},
\ref{hyp:morosowfork},} and~{\rm \ref{hyp:phimap}}  hold{\rm .}
If $\dorbit \in \dborbit(0)$ such that 
$$(X + \gg_{x,r^+}) \cap \dorbit \neq \emptyset,$$
then $\lsup{G}X \subset \overline{\dorbit}${\rm .}
\endproclaim

\vglue8pt

Here the closure is taken in the $p$-adic topology on $\gg$.

\demo{{P}roof {\rm (}\/J.\/{\rm -}\/L.~Waldspurger\/{\rm )}}
There exist $h \in G_x^+$ and $C \in C_{\gg_{x,r^+}}(Y)$ such that
$\lsup{h}(X + C) \in \dorbit$.  Thus, $X + C \in \dorbit$.  Since
$X+C$ is nilpotent, there exists (from Hypothesis~\ref{hyp:morosowfork})
a $\mu \in \X_*^k(\bG$) such that $\lsup{\mu(t)}(X+C) = t^2 \cdot (X+C)$ for
all $t \in k^\times$.  Since $C \in \oplus_{i \leq 0} \gg(i)$, we have
\begin{eqnarray*}
\lim_{t \rightarrow 0} \lsup{\lambda(t)\inv \mu(t)} (X + C) &=& X + \, \lim_{t
\rightarrow 0} \lsup{\lambda  (t)\inv} (t^2 \cdot C) \\
&=& X. \\ \noalign{\vskip-36pt}
\end{eqnarray*}
\phantom{soon}\hfill\qed
\enddemo

We close with a corollary to the above corollary.

\specialnumber{5.2.5}\proclaim{{C}orollary} \label{cor:uniqueorbitfa}
Assume that Hypotheses~{\rm \ref{hyp:decompnew},
\ref{hyp:morosowfork},} and~{\rm \ref{hyp:phimap}} hold{\rm .}  Choose $F^* \in
{\cal  F}(r)$ such 
that $F^* \subset \BB(Y,H,X)${\rm .}
If $\dorbit \in \dborbit(0)$ such that 
$$(X + \gg_{F^*}^+) \cap \dorbit \neq \emptyset,$$
then $\lsup{G}X \subset \overline{\dorbit}${\rm .}
\endproclaim

\demo{Proof}
Note that in this entire subsection the only assumption on $x$ was that
$x \in \BB(Y,H,X)$.  The result follows from
Corollary~\ref{cor:uniqueorbit}.
\enddemo

5.3. {\it A map for degenerate cosets}.  We assume that all of the hypotheses stated in~Section~4.2 hold.

The set ${I_r}$ is  too large.  We first restrict to degenerate cosets
of depth $r$.

\specialnumber{5.3.1}\numbereddemo{Definition}
\vglue2pt
\centerline{${\displaystyle I_r^n := \{(F^*,v) \in I_r \, | \, v \hbox{ is a degenerate element
of } V_{F^*}\}.}$} 
\enddemo

\specialnumber{5.3.2}\numbereddemo{{R}emark}
Suppose  $(F_i^*,v_i) \in I_r$ for $i = 1,2$ and  $(F_1^*,v_1) \sim
(F_2^*,v_2)$.  We have 
$(F_1^*,v_1) \in I_r^n$ if and only if $(F_2^*,v_2) \in I_r^n$.
\enddemo

Suppose that $(F^*,e) \in I_r^n$.  We wish to attach to $(F^*,e)$ a nilpotent
orbit $\dorbit(F^*,e) \in \dborbit(0)$. 

Suppose $x \in F^*$.  
We adopt the following conventions.  
If $e$ is trivial, then we declare that the ${\rmsl}_2(\ff)$-triple
$(f,h,e) \in 
V_{x,-r} \times V_{x,0} \times V_{x,r}$  completing $e$ is the trivial
triple.  Moreover, 
given a trivial ${\rmsl}_2(\ff)$-triple $(f,h,e)$ as above, we declare
that the ${\rmsl}_2(k)$-triple 
lifting $(f,h,e)$ is the trivial ${\rmsl}_2(k)$-triple. 

\specialnumber{5.3.3}\proclaim{Lemma} \label{lem:backoftt}
Suppose all the Hypotheses of Section~{\rm 4.2} hold{\rm .}  Suppose
$(F^*,e) \in I_r^n${\rm .}  
\begin{itemize}
\ritem{1.}
Fix $x \in F^*${\rm .}  There exists an ${\rmsl}_2(\ff)$\/{\rm -}\/triple
$(f,h,e) \in 
V_{x,-r} \times V_{x,0} \times V_{x,r}$ completing $e$ and an
${\rmsl}_2(k)$\/{\rm -}\/triple $(Y,H,X)$ which lifts $(f,h,e)${\rm .}
\ritem{2.}   
For any $x \in F^*${\rm ,} for any ${\rmsl}_2(\ff)$\/{\rm -}\/triple
$(f,h,e) \in 
V_{x,-r} \times V_{x,0} \times V_{x,r}$ completing $e${\rm ,} and for any
${\rmsl}_2(k)$\/{\rm -}\/triple $(Y,H,X)$ which lifts $(f,h,e)$ we have $F^*
\subset \BB(Y,H,X)$ and $\lsup{G}X$ is the unique nilpotent orbit of
minimal dimension which intersects the coset $e$
nontrivially{\rm . }
\end{itemize}

\endproclaim

{\it Proof}.
We first prove~(1).  Fix $x \in F^*$.  If $e$ is trivial, there is nothing to do.  Suppose
$e$ is nontrivial.
Hypothesis~\ref{hyp:pforallnew} says 
that there exist $f  \in V_{x,-r}$ and  $h \in V_{x,0}$
such that $(f,h,e)$
is an ${\rmsl}_2(\ff)$-triple (under the inherited Lie algebra
operation).
From Corollary~\ref{cor:liftsexist} we know that a
lift $(Y,H,X)$ of $(f,h,e)$ exists.

Now we prove~(2).  Suppose $x \in F^*$, the  ${\rmsl}_2(\ff)$-triple
$(f,h,e) \in 
V_{x,-r} \times V_{x,0} \times V_{x,r}$ 
completes $e$, and $(Y,H,X)$ is an 
${\rmsl}_2(k)$-triple which lifts $(f,h,e)$.
We have
$F^* \subset \BB(Y,H,X)$.
It follows from Corollary~\ref{cor:uniqueorbitfa} that $\lsup{G}X$ is the
unique nilpotent orbit 
of minimal dimension which intersects the coset $e$
nontrivially.  
\phantom{yoo}\hfill\qed
\vglue12pt

The following definition now makes sense.

\specialnumber{5.3.4}\numbereddemo{Definition} Suppose all the Hypotheses of Section~4.2 hold.
For $(F^*,e) \in I_r^n$ let $\dorbit(F^*,e)$ denote the
unique nilpotent orbit of minimal dimension which intersects the coset $e$
 nontrivially.
\enddemo

\specialnumber{5.3.5}\numbereddemo{{R}emark} \label{rem:mapginv}
If $g \in G$ and $(F^*,e) \in I_r^n$, then $\dorbit(gF^*,\lsup{g}e) =
\dorbit(F^*,e)$. 
\enddemo

5.4. {\it The map is well defined}.
Recall the equivalence relation on $I_r$ defined
in Section~3.6.  

\specialnumber{5.4.1}\proclaim{Lemma}  \label{lem:welldefined}
We assume that all of the hypotheses of~Section~{\rm 4.2} hold{\rm .}
The map from $I_r^n$ to $\dorbit(0)$ which sends $(F^*,e)$
to $\dorbit(F^*,e)$ induces a  well\/{\rm -}\/defined
map from $I_r^n / \!\sim$ to $\dorbit(0)${\rm .}
\endproclaim

\demo{Proof}
Suppose $(F_i^*,e_i) \in I_r^n$ for $i = 1,2$.  We need to show that
if $(F_1^*,e_1) \sim (F_2^*,e_2)$, then $\dorbit(F^*_1,e_1) =
\dorbit(F_2^*, e_2)$.  We may assume that $e_i \in V_{F_i^*}$ is not trivial.

Choose $x_i \in C(F_i^*)$.  Since $(F_1^*,e_1) \sim (F_2^*,e_2)$,
there exist $g \in G$ and an apartment $\AA$ in $\BB(G)$  such that
$$\emptyset \neq A(\AA, F_1^*) = A(\AA,gF_2^*)$$
and
$$e_1 \ident \lsup{g}e_2 \hbox{ \, in \,} V_{F_1^*} \ident V_{gF_2^*}.$$

From Remark~\ref{rem:mapginv} we can assume that $g = 1$.

Let $\bS$ denote the maximal $k$-split torus in $\bG$ corresponding to
$\AA$.  Let $\bfS$ denote the maximal $\ff$-split torus in $\bfG_{x_1}$
corresponding to $\bS$.

Complete $e_1$ to an  ${\rmsl}_2(\ff)$-triple  $(f_1,h_1,e_1)
\in V_{x_1,-r} \times  V_{x_1,0} \times V_{x_1,r}$ and suppose
$\bar{\lambda} \in \X_*^\ff(\bfG_{x_1})$ is adapted to this triple.  There
exists $h \in G_{x_1}$ such that $\lsup{\bar{h}}\bar{\lambda} \in
\X_*(\bfS)$.  (Here $\bar{h}$ denotes the image of $h$ in
$\bfG_{x_1}(\ff)$.)  Since $F_1^*$ and $F_2^*$ are strongly
$r$-associated and $e_1  \ident e_2$ in $V_{F_1^*} \ident V_{F_2^*}$,
it follows from Lemma~\ref{lem:identaut} that there exists $h' \in
G_{x_1} \cap G_{x_2}$ such that
$$\lsup{h}e_1 \ident \lsup{h'}e_1 \ident \lsup{h'}e_2 \hbox{\, in \,}
V_{F_1^*} \ident V_{F_1^*} \ident V_{F_2^*}.$$
Let $\lambda \in \X_*(\bS)$ be the lift of
$\lsup{\bar{h}}\bar{\lambda}$.  

Let $\gg = \mathbold{\oplus}_j \gg(j)$ be the decomposition of $\gg$ arising from
$\lambda$.  We have $\gg_{F_i^*} = \mathbold{\oplus}_j \gg_{F_i^*}(j)$.  There exists
an $X \in \gg_{F_1^*}(2) \cap \gg_{F_2^*}(2)$ such that the image of $X$
in $V_{F_i^*}$ is $\lsup{h'}e_i$.

It follows from Corollary~\ref{cor:liftsexist} and
Lemma~\ref{lem:backoftt} that
\vglue4pt \hfill
${\displaystyle \dorbit(F_i^*,e_i) = \dorbit(F_i^*,\lsup{h'}e_i) = \lsup{G}{X}.} $
 \enddemo
\vglue6pt

5.5. {\it Distinguished cosets}. 
We assume that all of the hypotheses of~Section~4.2 hold.

The set  $I_r^n/ \!\sim$ is too large.  We now restrict our attention to
{\it distinguished} cosets of depth $r$.

\specialnumber{5.5.1}\numbereddemo{Definition}\label{def:distinguished}
We define $I_r^d \subset I_r^n$ to be those pairs $(F^*,e) \in I_r^n$
such that 
for any $x \in F^*$, for any
${\rmsl}_2(\ff)$-triple 
$(f,h,e) \in V_{x,-r} \times V_{x,0} \times V_{x,r}$ completing $e$,
and for any ${\rmsl}_2(k)$-triple $(Y,H,X)$ which lifts $(f,h,e)$,
we have
that $F^*$ is a maximal generalized $r$-facet in $\BB(Y,H,X)$.
\enddemo

\specialnumber{5.5.2}\numbereddemo{{R}emark}  If $r = 0$, then it can be shown that this definition of
distinguished is equivalent to the usual one.  That is, if $(F^*,e)
\in I_0^d$, then $e$ does not lie in any proper Levi subalgebra of the
$\ff$-{L}ie algebra $V_{F^*,0}$.
\enddemo

\specialnumber{5.5.3}\proclaim{Lemma}  
Suppose all of the hypotheses of Section~{\rm 4.2} hold{\rm .}
If $(F^*,e) \in I_r^n$ and $e$ is nontrivial{\rm ,} then  $(F^*,e) \in I_r^d$
if and only if there exist
an $x \in F^*${\rm ,}
an ${\rmsl}_2(\ff)$\/{\rm -}\/triple $(f,h,e) \in V_{x,-r} \times
V_{x,0} \times V_{x,r}$ completing $e${\rm ,} and
an ${\rmsl}_2(k)$\/{\rm -}\/triple $(Y,H,X)$ in $\gg$
lifting $(f,h,e)$
such that $F^*$ is a maximal generalized $r$\/{\rm -}\/facet in $\BB(Y,H,X)${\rm .} 
\endproclaim

\demo{Proof}
``$\Rightarrow$'':  This follows from the definitions.

``$\Leftarrow$'':  Suppose we have
an $x \in F^*$,
an ${\rmsl}_2(\ff)$-triple $(f,h,e) \in V_{x,-r} \times
V_{x,0} \times V_{x,r}$ completing $e$, and
an ${\rmsl}_2(k)$-triple $(Y,H,X)$ in $\gg$
lifting $(f,h,e)$
such that $F^*$ is a maximal generalized $r$-facet in $\BB(Y,H,X)$.

Suppose we also have data 
$x' \in F^*$, 
an ${\rmsl}_2(\ff)$-triple 
$(f',h',e) \in V_{x',-r} \times V_{x',0} \times V_{x',r}$ completing $e$,
and an ${\rmsl}_2(k)$-triple $(Y',H',X')$ 
in $\gg$ which lifts $(f',h',e)$
such that $F^*$ is not a maximal generalized $r$-facet in $\BB(Y',H',X')$.
We will derive a contradiction.

Since $\gg_{x',-r} = \gg_{x,-r}$ (from Lemma~\ref{lem:needed}) and
$\gg_{x',r} = \gg_{x,r}$, we have  
$$[\gg_{x',-r},\gg_{x',r}] \subset \gg_{x,0}.$$  Thus, we can
assume that $x = x'$.

From Lemma~\ref{lem:backoftt} we have
$$\lsup{G}X = \lsup{G}X' = \dorbit(F^*,e).$$
From Corollary~5.2.3 we can assume (after replacing
$(Y',H',X')$ with a $G_{x}^+$-conjugate) that $X = X'$.  From
Hypothesis~\ref{hyp:morosowfork}, 
there exists a $g
\in C_{G}(X)$ such that $Y' = \lsup{g}Y$ and $H' = \lsup{g}H$.
Consequently, $\BB(Y',H',X) = g \BB(Y,H,X)$.

By assumption, $g\inv F^* \subset g\inv  \BB(Y',H',X) = \BB(Y,H,X)$ is
not a maximal 
generalized $r$-facet in $\BB(Y,H,X)$.    Since $\dim g\inv F^* =
\dim F^*$, this is a contradiction.
\enddemo

\specialnumber{5.5.4}\numbereddemo{{R}emark}
Suppose $(F_i^*,e_i) \in I_r^n$ for $i = 1,2$ and $(F_1^*,e_1) \sim
(F_2^*,e_2)$.  From Lemma~\ref{lem:welldefined} we have
$\dorbit(F^*_1,e_1) = \dorbit(F^*_2,e_2)$.  From the
proof above, we conclude that $(F_1^*,e_1) \in I_r^d$ if and only if
$(F_2^*,e_2) \in I_r^d$.
\enddemo

5.6. {\it A bijective correspondence}.
We assume that all of the hypotheses of~Section~4.2 hold.
In this subsection we establish the following theorem.

\specialnumber{5.6.1}\proclaim{Theorem} \label{thm:parametrize}  
Assume that all of the hypotheses of Section~{\rm 4.2} hold{\rm .}
There is a bijective correspondence between $I_r^d / \!\sim$ and
$\dorbit(0)$ given by the map which sends $(F^*,e)$ to $\dorbit(F^*,e)${\rm .}
\endproclaim

{\it Proof}.
We have already seen that the map is well defined.  

We first show that
the map is injective.
We need to show that if  $(F_1^*,e_1),\break (F^*_2,e_2) \in I_r^d$ and
$\dorbit(F_1^*,e_1) = \dorbit(F_2^*,e_2)$, then $(F_1^*,e_1) \sim
(F_2^*,e_2)$.

If $\dorbit(F_i^*,e_i) = \{0\}$, then $F_i^*$ is open in $\BB(G)$ and
the result 
follows.  Thus, we may assume that $e_i$ is not trivial.

Fix $x_i \in C(F_i^*)$.  Complete $e_i$ to an ${\rmsl}_2(\ff)$-triple
$(f_i,h_i,e_i) \in V_{x_i,-r}\break \times V_{x_i,0} \times V_{x_i,r}$.
From Corollary~\ref{cor:liftsexist} we may lift $(f_i,h_i,e_i)$ to an
${\rmsl}_2(k)$-triple 
$(Y_i,H_i,X_i)$ in $\gg$.   From Lemma~\ref{lem:backoftt} we have
$\dorbit(F_i^*,e_i) = 
\lsup{G}X_i$.  

Since $\dorbit(F_1^*,e_1) = \dorbit(F_2^*,e_2)$, there exists a $g \in
G$ such that $(\lsup{g}Y_2, \lsup{g}H_2, \lsup{g}X_2)\break =
(Y_1,H_1,X_1)$.  Consequently, since $(F_i^*,e_i) \in I_r^d$, from
Lemma~\ref{lem:maximalgenr} we have that $F_1^*$ and $gF_2^*$ are
strongly $r$-associate.  Thus, there exists an apartment $\AA$ in
$\BB(G)$ such that
$$\emptyset \neq A(\AA,F_1^*) = A(\AA,gF_2^*).$$
Moreover, since $X_1$ has image $e_1$ in $V_{F_1^*}$ and $X_1$ has image
$\lsup{g}e_2$ in $V_{gF_2^*}$, we have $X_1 \in \gg_{F^*_1} \cap
\gg_{gF^*_2}$ and 
$$e_1 \ident \lsup{g}e_2 \hbox{\, in \,} V_{F_1^*} \ident V_{gF_2^*}.$$
Thus, the map is injective.

We now show that the map is surjective.
Suppose $\dorbit \in \dborbit(0)$.  

If $\dorbit$ is trivial, then let
$F^*$ be  an open generalized $r$-facet and let $e$ be trivial in
$V_{F^*}$. We have $(F^*,e) \in I_r^d$ and \pagebreak $\dorbit(F^*,e) = \{0\}$.

Suppose $\dorbit$ is not trivial.
Fix $X \in \dorbit$.  Complete $X$
to an ${\rmsl}_2(k)$-triple 
$(Y,H,X)$ in $\gg$.  Let $F^*$ be a maximal generalized $r$-facet in
$\BB(Y,H,X)$ and let $e$ denote the image of $X$ in $V_{F^*}$.  We
will be done if we can show that $\dorbit(F^*,e) = \lsup{G}X$.  This,
however, follows from Lemma~\ref{lem:backoftt}~(2).
\hfill\qed\vglue9pt

For future reference, we record the following corollary of the proof of
Theorem~\ref{thm:parametrize}.
\vglue6pt

{\elevensc Corollary 5.6.2.}
{\it Assume that all of the hypotheses of Section~{\rm 4.2} hold{\rm .}
Suppose $(F_1^*,e_1), (F_2^*,e_2) \in I_r^d$ and $(F_1^*,e_1) \sim  
(F_2^*,e_2)${\rm .}   There exists $g \in G$ and an
${\rmsl}_2(k)$\/{\rm -}\/triple $(Y,H,X)$ such that
\vglue4pt
{\rm 1.} $X \in \gg_{F^*_1} \cap \gg_{gF_2^*}${\rm ,}
\vglue4pt
{\rm 2.} $X$ has image $e_1$ in $V_{F_1^*}$ and image $\lsup{g}e_2 \in
V_{gF_2^*}${\rm ,} and
\vglue4pt
{\rm 3.} $F_1^*$ and $gF_2^*$ are maximal generalized $r$-facets in
$\BB(Y,H,X)${\rm .}  }
 \vglue12pt

  \centerline{\bf Appendix A. Some comments on Hypothesis~\ref{hyp:pforallnew}}

\vglue4pt {\rm A.1.} {\it Introduction}.
We will show that, subject to some conditions on $k$ and $\bG$,
Hypothesis~\ref{hyp:pforallnew} is valid.  This result is 
related to material  
in~\cite{howe-moy:minimal} and\break \cite[\S2]{morris:fundamental} . No
attempt has been 
made to produce an optimal set of conditions.

Fix $x$, $r$, and $X$ as in the statement of
Hypothesis~\ref{hyp:pforallnew}.  Without loss of generality, we
assume throughout this 
appendix that $\bG$ is connected.
\vglue4pt

A.2. {\it An ${\rmsl}_2(\ff)$-triple}.  In this subsection, we establish
the existence of $Y$ and $H$ in $\gg$ satisfying the requirements of 
Hypothesis~\ref{hyp:pforallnew}.

 We let $\gg'$ (resp., $\zz$) denote
the Lie algebra 
of the group of $k$-rational points of the the derived group of $\bG$
(resp., the connected component of the
center of $\bG$).   From~\cite[Proposition~3.1]{adler-roche:intertwining}
if $p$ is larger than some constant which may be determined by
examining the absolute root datum of $\bG$, then we may
assume that $\gg = \zz + 
\gg'$.   Thus, without loss of
generality, we may assume that $\bG$ 
is semisimple.

We claim that, 
under suitable conditions on $k$ and $\bG$, the Killing form $\kappa$
identifies 
$\gg_{x,s}$ with 
$\gg^*_{x,s}$ 
for all $s \in \R$; in particular,
for all $Z \in \gg_{x,s}
\smallsetminus \gg_{x,s^+}$ there exists a $W \in \gg_{x,-s}
\smallsetminus \gg_{x,(-s)^+}$ such that $\kappa(Z,W) \in R^\times$. 
Indeed,   
since $\gg_{x,s}$\break \cite[Proposition~1.4.1]{adler:thesis} and $\kappa$
behave well with respect to Galois 
descent, we may reduce to the case when $\bG$ is $k$-split.  If $\bG$ is
$k$-split, we can fix a  
Chevalley basis for $\rmg$.  In this situation the statement follows
if $p$ is greater than some constant which may be derived from the
absolute root datum of $\bG$.    

\pagegoal=49pc
Recall the definition of the  finite-dimensional $\ff$-Lie algebra
$\overline{\gg}_x$ from Section~4.2.  From the
previous paragraph we see that, under suitable 
conditions on $k$ and $\bG$, the representation 
${\dad}$ 
of $\overline{\gg}_x$ has a nondegenerate trace-form. Let $e \in
\overline{\gg}_x$ denote the image of $X$ in
$V_{x,r}$.
From Hypothesis~\ref{hyp:decompnew} we have that
${\dad}(e)^{m} = 0$ for some $m \leq (p-2)$.
Thus from~\cite[Proposition~5.3.1]{carter:finite} there exists an 
${\rmsl}_2(\ff)$-triple $(f,h,e)$ in $\overline{\gg}_x$ completing $e$.\pagebreak

We claim that we may assume that $f \in V_{x,-r}$ and $h
\in V_{x,0}$.  (We already know that $e \in V_{x,r}$.)    Indeed, let $f_{-r}$
denote the image of $f$ under the 
projection of $\overline{\gg}_x$ onto
$V_{x,-r} \subset \overline{\gg}_x$.  Since $e 
\in V_{x,r} \subset \overline{\gg}_x$, it follows that
$\ad(e)^2 f_{-r} = \ad(e)^2 f = -2e$.  Let $h_0$ denote the image of $h$
under the projection of $\overline{\gg}_x$ onto $V_{x,0}$.  We have
$\ad(e)f_{-r} = h_0$ and $\ad(e) h_0 = -2e$.  
From Morosov's 
theorem (see~\cite[Lemma~7, p.~98]{jacobson:lie}
or~\cite[the proof of Proposition~5.3.1; in particular, 
pp.~140--141]{carter:finite}) there exists $f' \in
\overline{\gg}_x$ such that  
$(f',h_0,e)$ is an ${\rmsl}_2(\ff)$-triple.  Since $h_0 \in V_{x,0}$,
$e \in V_{x,r}$, and $[V_{x,s}, V_{x,t}] \subset V_{x,(s+t)}$ for all
$s,t \in \R$, we can assume
that $f' \in V_{x,-r}$.  Thus we may assume that $f \in V_{x,-r}$ and $h
\in V_{x,0}$.  
\pagegoal=48pc

Consequently, we can choose lifts
$Y \in \gg_{x,-r}$ of $f$ and $H \in \gg_{x,0}$ of $h$ which
satisfy the initial requirements of Hypothesis~\ref{hyp:pforallnew}.  
From Lemma~\ref{lem:degbenil} we may assume that $Y$ is nilpotent.
Thus from Hypothesis~\ref{hyp:decompnew} we have that
${\dad}(f)^{m} = {\dad}(e)^{m} = 0$ for some $m \leq
(p-2)$.
Therefore, from~\cite[Proposition~5.4.8]{carter:finite} we have that
$\overline{\gg}_x$ is a direct sum of irreducible $(f,h,e)$-modules of
highest weight at most $(p-3)$.

\demo{{\rm A.3.} A one-parameter subgroup in $\bfG_x$}
We now wish to establish the existence of $\bar{\lambda} \in
\X_*^\ff(\bfG_x)$ satisfying the requirements of
Hypothesis~\ref{hyp:pforallnew}.  We continue to
use the notation introduced above, but  
we now remove the assumption that $\bG$ is semisimple.

Since $h$ is semisimple, from~\cite[Proposition~11.8 and its
proof]{borel:linear} there exists a maximal $\ff$-torus $\bfT$ of
$\bfG_x$ such that $h \in \Lie(\bfT)$. 
Let $\bT$ be a maximal $K$-split $k$-torus of $\bG$ associated to
$\bfT$. (That is, we have $x \in \AA(\bT,K)$ and the image of
$\bT(K) \cap \bG(K)_x$ in $\bfG_x(\ffc)$ is $\bfT(\ffc)$.  That the
torus $\bT$ exists follows from the argument made in the final paragraph of
the proof of~\cite[Proposition~5.1.10]{bruhat-tits:two}.)  

Let
$$\overline{\gg}(\ffc)_x := \bigoplus_{s \in \R/ {\ell \cdot \Z}} {\rmg }(K)_{x,s}/{\rmg }(K)_{x,s^+}.$$
As in~Section~4.2, we give $\overline{\gg}(\ffc)_x$ a natural Lie
algebra structure; in fact, $\overline{\gg}_x$ may be identified with
the set of 
$\Gal(\ffc/\ff)$-fixed points of the $\ff$-Lie algebra
$\overline{\gg}(\ffc)_x$.
Since $\overline{\gg}(\ffc)_x$ decomposes into irreducible
$(f,h,e)$-modules, it follows from~\cite[Lemmas~5.5.3 and
 5.5.4]{carter:finite} that there exists a one-parameter subgroup
$\bar{\lambda}_2 \colon \glo \rightarrow {\rm GL}(\overline{\gg}(\ffc)_x)$
defined over 
$\ff$ such that  for $v \in \overline{\gg}(\ffc)_x$
$$\hbox{if ${\bar{\lambda}_2(t)} \cdot v = t^iv$, then $\dabs{i} \leq
(p-3)$ and 
$[h,v] = iv$.}$$

For
$\alpha \in \Phi(\bT,K)$, we denote by $(\overline{\gg}(\ffc)_x)_\alpha$ the
(nontrivial) subspace of $\overline{\gg}(\ffc)_x$ on which
$\bfT$ acts by $\alpha$.  
We define a linear map $\lambda_2$
from the $K$-root lattice in $\X^*(\bT)$ to
$\Z$ via $\bar{\lambda}_2$.  That is, for $\alpha \in \Phi(\bT,K)$ we
define $\langle \lambda_2, \alpha 
\rangle \in \Z$ by the equality
$${\bar{\lambda}_2(t)} \cdot v = t^{\langle \lambda_2, \alpha
\rangle}v$$
for all $v \in (\overline{\gg}(\ffc)_x)_\alpha$ and extend linearly.
Note that $\lambda_2$ is $\Gal(K/k)$-invariant.

For $i \in \Z$, we define
$${\rmg }(K)(i) := \bigoplus_{\alpha \in \Phi(\bT,K) \cup \{0\}; \langle
\lambda_2 , \alpha 
\rangle = i} {\rmg }(K)_\alpha.$$
Write $\gg(i)$ for the $\Gal(K/k)$-fixed points of ${\rmg }(K)(i)$.  Let
$X \in \gg(2) \cap \gg_{x,r}$ be any lift of $e$.  Since
$\bar{\gg}_x$ decomposes into irreducible $(e,h,f)$-modules of
highest weight at  most $(p-3)$, we conclude that $\dad(e)^2 \colon
\bar{\gg}_x(-2) \rightarrow \bar{\gg}_x(2)$ is an  isomorphism (here
$\bar{\gg}_x(\pm2) = \{v \in \bar{\gg}_x \, | \, \dad(h)v = \pm2v
\}$).  Thus, since $k$ is complete, for all $s \in \R$ the map
$\dad(X)^2 \colon \gg(-2) \cap \gg_{x,s-r} \rightarrow \gg(2) \cap
\gg_{x,s+r}$ is an isomorphism (see also the proof of
Lemma~\ref{lem:Xsquareiso}).  Hence, there exists $Y \in \gg(2) \cap
\gg_{x,-r}$ such that $\dad(X)^2Y = -2X$.  Let $H = \dad(X)Y$.  Note
that $H$ is a lift of $h$ and $Y$ is a lift of $f$.  We claim that
$(Y,H,X)$ is an ${\rmsl}_2(k)$-triple.  From Hypothesis~\ref{hyp:decompnew}
and Morosov's theorem, there exists $Y' \in \gg$ such that
$(Y',H,X)$ is an ${\rmsl}_2(k)$-triple.  However, as usual, we conclude
that we may assume that $Y' = Y$ (see also the proof of
Corollary~\ref{cor:liftsexist}). 

From Hypothesis~\ref{hyp:morosowfork} there exists a $k$-homomorphism
$\varphi \colon {\bf SL}_2 \rightarrow \bG$  such that  $d\varphi {\big(
\begin{array}{cc}\scrs{0}\hskip-5pt&\scrs{1}\\[-6pt]\scrs{0}\hskip-5pt&\scrs{0}\end{array} \big)} = X$, $d\varphi
{\big(
\begin{array}{cc}\scrs{0}\hskip-5pt&\scrs{0}\\[-6pt]\scrs{1}\hskip-5pt&\scrs{0}\end{array} \big)} = Y$, and $d\varphi
{\big(
\begin{array}{cc}\scrs{1}\hskip-5pt&\scrs{0}\\[-6pt]\scrs{0}\hskip-5pt&\scrs{-1}\end{array} \big)} = H$.  Let
$\lambda$ denote the one-parameter subgroup $t \mapsto
\varphi{\big(
\begin{array}{cc}\scrs{t}\hskip-5pt&\scrs{0}\\[-6pt]\scrs{0}\hskip-5pt&\scrs{t^{-1}}\end{array} \big)}$.    Statements
(1) and (2) of Hypothesis~\ref{hyp:morosowfork} along with
Hypothesis~\ref{hyp:decompnew}  imply that for all $v \in \gg$
\begin{equation} 
\speqnu{*}
\label{equ:star}
\hbox{if $\lsup{\lambda(t)}v = t^i v$, then $\dabs{i} \leq (p-3)$ and
$\dad(H)v = iv$}.
\end{equation}
 
From Corollary~\ref{cor:liecentisfixm} we have $x \in
\BB(C_G(\lambda))$.  Thus, there exists a maximal $k$-split torus
$\bS$ such that $x \in \AA(\bS,k)$ and $\lambda \in \X_*(\bS)$.
Let $\bar{\lambda}$ denote the image of ${\lambda}$ in
$\X_*^\ff(\bfG_x)$.  The image of $d\bar{\lambda}$ in $\Lie(\bfG_x)$
coincides with the one-dimensional subspace spanned by $h$.   
From
(\ref{equ:star}) and that fact that $H$ is a lift of $h$ we have that
for all $v \in \bar{\gg}_x$
$$\hbox{if $\lsup{\bar{\lambda}(t)}v = t^i v$, then $\dabs{i} \leq (p-3)$ and
$\dad(h)v = iv$}.$$

Finally, we consider the uniqueness statement.  
Fix $i \in \Z$ such that $-2 \leq i \leq 2$.  Note that if $m \in \N$
and $v \in \overline{\gg}_x$ such that $\dad(h)v =  iv$, then
$\lsup{\bar{\lambda}(t)}\bigl(\dad(e)^mv\bigr) = t^{i+2m} \bigl(\dad(e)^mv\bigr)$ and
$\lsup{\bar{\lambda}(t)}\bigl(\dad(f)^mv\bigr) = t^{i -2m}
\bigl(\dad(f)^mv\bigr)$.  Since $\overline{\gg}_x$ is spanned by the set of all
vectors of the form 
$\dad(e)^mv$ or $\dad(f)^mv$ ($m$ and $v$ as above), we conclude that
$\bar{\lambda}$ is uniquely determined up to an element of
$\X_*(\bfZ_x)$ whose differential is zero.

\end{document}